# SPATIAL HOMOGENIZATION IN A STOCHASTIC NETWORK WITH MOBILITY

By Florian Simatos and Danielle Tibi[1]

*INRIA and Université Paris 7*

A stochastic model for a mobile network is studied. Users enter the network, and then perform independent Markovian routes between nodes where they receive service according to the Processor-Sharing policy. Once their service requirement is satisfied, they leave the system. The stability region is identified via a fluid limit approach, and strongly relies on a "spatial homogenization" property: at the fluid level, customers are instantaneously distributed across the network according to the stationary distribution of their Markovian dynamics and stay distributed as such as long as the network is not empty. In the unstable regime, spatial homogenization almost surely holds asymptotically as time goes to infinity (on the normal scale), telling how the system fills up. One of the technical achievements of the paper is the construction of a family of martingales associated to the multidimensional process of interest, which makes it possible to get crucial estimates for certain exit times.

**1. Introduction.** Recent wireless technologies have triggered interest in a new class of stochastic networks, called *mobile networks* in the technical literature [3, 11]. In contrast with Jackson networks where users move upon completion of service at some node, in these mobile networks, transitions of customers within the network occur *independently* of the service received. Moreover, at any given time, each node capacity is divided between the users present, whose service rate thus depends on the capacity and on the state of occupancy of the node. Once his initial service requirement has been fulfilled, a customer definitively leaves the network. In [3], complex capacity sharing policies are considered, but in the simplest setting, which will be of interest to us, nodes implement the Processor-Sharing discipline by dividing their capacity equally between all the users present. Previous works [3, 11] have mainly focused on determining the stability region of such networks, and it has been commonly observed that the users' mobility represents an opportunity for the network to increase this region. Indeed, because of their mobility, users offer a diversity of channel conditions to the base stations (in charge of allocating the resources of the nodes), thus allowing them to select the users in the most favorable state. Such a scheduling strategy is sometimes referred to as an *opportunistic scheduling* strategy (see [2] and the references therein for more details).

[1]On leave at INRIA Paris-Rocquencourt.

In the present paper, we investigate from a mathematical standpoint a basic Markovian model for a mobile network derived from [3]. In this simple setting, customers arrive in the network according to a Poisson process with intensity $\lambda$, and move independently within the network, according to some Markovian dynamics with a common rate matrix $Q$. Service requirements are exponentially distributed with mean 1, and customers are served at each node they visit according to the Processor-Sharing discipline, until their demand has been satisfied. The total capacity of the network, defined as the sum of all the individual capacities of the nodes, is denoted by $\mu$. It corresponds to the instantaneous output rate of the network when no node is empty, that is, when there is at least one customer at each node.

It is of particular interest to note that, even if $Q$ is reversible, because of the arrival and departure processes, the system is not reversible. This contrasts with earlier works in which particle systems with similar dynamics have been investigated under reversibility assumptions. In [6], the authors look at a closed system (i.e., with parameters $\lambda = \mu = 0$) where transition rates are chosen such as to yield a reversible dynamics. In this case, the stationary distribution of the system has a product form, and the authors are interested in showing that the convergence to equilibrium is exponentially fast. Their approach essentially relies on logarithmic Sobolev type inequalities.

In our case, however, a different set of questions is addressed, involving different tools. Since the system under consideration is open, it may be unstable, so that a natural issue is to determine the stability region. We prove, as was conjectured in [3], that the intuitive, simple condition $\lambda < \mu$ is indeed the stability condition (the critical case $\lambda = \mu$ is not considered). In contrast with Jackson networks for which the stability condition is local, in the sense that each node has to satisfy some constraint, here only the global quantities $\lambda$ and $\mu$ matter. This shows that mobility allows to make the most of the potential service capacity of the network, corroborating the results previously mentioned. Note that $\lambda < \mu$ being a necessary condition is obvious, since $\mu$ is the maximal output rate. But surprisingly, proving that it is sufficient requires very technical tools, including the use of fluid limits and martingale techniques. In particular, the long and tedious Appendix is solely devoted to the construction of a martingale which provides key estimates for showing that $\lambda < \mu$ corresponds to a stable system.

This martingale is a multidimensional (therefore complicated) generalization of the martingale built in [10] for the $M/M/\infty$ queue, and this is not completely surprising, since as will be seen, the model inherits salient properties of the $M/M/\infty$ queue. Besides, the construction of a martingale associated to a multidimensional process represents one of the technical achievements of this paper: such examples are indeed pretty scarce in the literature. Similar to [10], the approach relies on building a family of space–time harmonic functions indexed by some parameter $c \in \mathbb{R}^n$, and then on integrating over $c$ in such a way as to preserve the harmonic property.

Through studying both the stability region and the unstable regime, a detailed description of the behavior of the system is given, resulting in two versions (stable and unstable) of the following rough property: when many users are present in the network, they get approximately distributed among the nodes according to the unique invariant distribution $\pi$ associated with $Q$, the latter being assumed irreducible. It must be emphasized that yet, contrary to [3], customers' movements are not assumed stationary.

As a first argument for this spatial homogenization, the law of large numbers suggests that, when the total number of users initially present in the network is large, the proportions of users at the different nodes should be close to $\pi$ after some time, related to the convergence to $\pi$ of the Markov process associated to $Q$. The more delicate question that next arises, of how long these proportions stay close to $\pi$, constitutes the main challenging issue of the paper, that requires martingale techniques for estimating the deviation time from $\pi$.

The short-term reach of $\pi$ is understandable from an analogy with the $M/M/\infty$ queue; indeed, independence of the customers' trajectories yields that, similar to the $M/M/\infty$ queue, the output rate from any node due to inner transitions is directly proportional, through $Q$, to the number of customers at this node. When the network is overloaded, the relative occupancies of the nodes should then, after a while, be close to the internal traffic balance ratios given by $\pi$.

A more explicit analogy with another classical queueing model is provided by the following simple but crucial observation: as long as no node is empty, the total number of customers simply evolves as an $M/M/1$ queue with input rate $\lambda$ and output rate $\mu$. This is the case, in particular, when the distribution of customers is close to $\pi$. This interplay between, on one hand, the proportions of customers at the different nodes, and, on the other hand, their total number, will underly the analysis throughout the paper.

While the short-term behavior, which results in the spreading of customers according to $\pi$, is dominated by the $M/M/\infty$ dynamics, the long-term behavior is essentially driven by the $M/M/1$ dynamics of the total number of customers. This naturally suggests that two different scalings have to be considered: one, corresponding to the $M/M/\infty$ dynamics where only space is scaled and not time; and a second one where both space and time are scaled corresponding to the fluid scaling of the $M/M/1$ queue. Note that the natural scaling for the $M/M/\infty$ queue is the so-called Kelly scaling, in which space and input rate are scaled. Here, since the input rate at each node due to inner transitions is a linear function of the numbers of customers at the different nodes, there is no need to scale the external input rate $\lambda$. Inner movements dominate the dynamics and the space scaled process converges, analogously to the $M/M/\infty$ queue under Kelly's scaling, to some deterministic trajectory, with limit point at infinity here given by $\pi$.

The coexistence of these two different scalings makes the use of fluid limits both original and challenging. Fluid limits are a standard tool in the analysis of complicated stochastic networks. Rybko and Stolyar [13] is one of the first papers that uses this technique together with Dai [8]. Dupuis and Williams [9] present similar ideas in the context of diffusions. In a series of papers, Bramson [4, 5] describes the precise evolution of fluid limits for various queueing networks. See also the books by Chen and Yao [7] and Robert [12]. In the context of networks, fluid limits have been used mainly for Markov processes which behave locally as random walks. For this reason, results related to fluid limits are sometimes presented as functional laws of large numbers. Because of the mixture of two different dynamics, given by the $M/M/1$ and $M/M/\infty$ models, our framework is somewhat different. A second important difference with the existing literature concerns tightness results which are usually easy to obtain, mainly because transition rates are generally bounded; this not the case here.

The long-term analysis is twofold. Deriving fluid limits requires a control on the process over time periods of the same order as the initial number of customers (since the fluid scaling parameter is the same for time and space). In the stable case this is obtained by showing that the deviation time from $\pi$ is essentially larger than the time for the underlying $M/M/1$ queue to empty. The unstable case exhibits a more striking behavior: the deviation time from $\pi$ is not only large compared to the initial number of customers, but is even *infinite* with high probability. This amounts to a control of the whole trajectory; the distribution of users among nodes stays trapped in any neighborhood of $\pi$ with high probability as the initial state is large. This result is related to a strong convergence result stating that, for any fixed (nonscaled) initial state, the system almost surely diverges along the direction of $\pi$. Note that a similar phenomenon has been exhibited in [1], in the context of branching Markov chains, that is, Galton–Watson branching processes where individuals located at some countable set of sites move at their birth time.

These various remarks and outline of results lead to the following organization for the paper. Section 2 gives a precise description of the stochastic model and introduces the notation that will hold throughout the paper. We have already mentioned the construction of a martingale which gives important estimates through optional stopping techniques. Section 3 introduces this martingale, and provides the main estimate that will be used. Due to its technicality, the construction of the martingale is postponed to the Appendix.

Section 4 establishes a decomposition of the process as, mainly, the difference between two processes of the same type but with no departures. For such a process (with null service capacity), a representation involving labelled particles is given. Both representations will help derive the almost sure convergence result of Section 6.

The three last sections are devoted to analyzing the behavior of the system. Section 5 deals with the short-term behavior, thus studying the only space renormalized process. Section 6 studies the supercritical case $\lambda > \mu$, establishing among other results the almost sure convergence of the proportions to the equilibrium distribution $\pi$ as $t \to \infty$. Finally, Section 7 proves the stability of the system in the subcritical case $\lambda < \mu$.

**2. Framework and notation.** This section gives a precise description of the model under consideration and introduces the main notation. The network is described by a Markov process $X = (X(t), t \geq 0)$ characterized by its infinitesimal generator, given by (1) below.

Section 6 will make use, in the particular case of null service capacity, of a more explicit representation of $X$ involving a sequence of Markov jump processes that represent the trajectories of the successive customers entering the network. The general description of the system through its Markovian dynamics provided in the present section is, however, sufficient for most results of the paper, especially for building a family of martingales and for determining the stability condition.

The network consists of $n$ nodes between which customers perform independent (continuous-time) Markovian routes during their service. In this setting, transitions of customers from one node to another are driven by some rate matrix $Q = (q_{ij}, 1 \leq i, j \leq n)$ and are thus not triggered by service completion.

New customers arrive at node $i = 1, \dots, n$ according to a Poisson process with intensity $\lambda_i \geq 0$, and then move independently according to the Markovian dynamics defined by $Q$. The arrival processes at the different nodes are independent, so that the global arrival process is Poisson with intensity $\lambda = \sum_1^n \lambda_i$. The case $\lambda = 0$ corresponds to a system with only initial customers, and no new arrivals.

Upon arrival, or at time $t = 0$ for those initially present, customers generate a service requirement which is exponentially distributed with mean 1. All service requirements, arrival processes and Markovian routes are assumed to be mutually independent.

Node $i$, $1 \leq i \leq n$, has service capacity $\mu_i \geq 0$, which is divided at any time between the customers present, according to the Processor-Sharing discipline: if $N$ is the number of customers present at node $i$, then each of these $N$ customers is served at rate $\mu_i/N$. The service rate of a given customer thus evolves in time, depending on his current position and on its occupancy level. Once a customer has received a service that meets his initial requirement, he leaves the network.

The total service capacity of the network is defined as $\mu = \sum_1^n \mu_i$. Notice that, due to the exponential nature of the services, the mechanism of departure from one node by completion of service does not distinguish the present Processor-Sharing discipline from the FIFO discipline: the instantaneous output rate from the system at node $i$ is $\mu_i$, provided that node $i$ is not vacant. The total output rate is then $\mu$ when no node is empty.

The process of interest is $X = (X(t), t \geq 0)$ defined by

$$X(t) = (X_1(t), \ldots, X_n(t)), \qquad t \geq 0,$$

where $X_i(t)$, for $i = 1, \ldots, n$, is the number of customers present at node $i$ at time $t$. The Markovian nature of the movements together with the exponential assumption for the service distribution imply that $X$ is a Markov process in $\mathbb{N}^n$ with infinitesimal generator $\Omega$ given, for any function $f : \mathbb{N}^n \to \mathbb{R}$ and any $x = (x_1, \ldots, x_n) \in \mathbb{N}^n$, by

$$\begin{aligned} \Omega(f)(x) = {} & \sum_{i=1}^n \lambda_i (f(x+e_i) - f(x)) \\ & + \sum_{i=1}^n \mathbb{1}_{\{x_i>0\}} \mu_i (f(x-e_i) - f(x)) \\ & + \sum_{1 \leq i \neq j \leq n} q_{ij} x_i (f(x+e_j-e_i) - f(x)), \end{aligned} \tag{1}$$

where $e_i \in \mathbb{N}^n$ has all coordinates equal to 0, except for the $i$th one, equal to 1.

The Introduction has highlighted that this system is a mixture of two classical models in queueing theory, the $M/M/1$ and the $M/M/\infty$ queues. This is readable in the expression of the generator given in (1) where the two first sums are reminiscent of the $M/M/1$ queue and the last one of the $M/M/\infty$ queue.

The rate matrix $Q$ is assumed to be irreducible, admitting $\pi = (\pi_i, 1 \le i \le n)$ as its unique stationary distribution characterized by the relation

$$\pi Q = 0.$$

For technical reasons related to the construction of the martingale introduced in Section 3 (see in the Appendix), we require the additional assumption that $Q$ is diagonalizable. This assumption is satisfied if $Q$ is reversible with respect to $\pi$, but it is in general a much less restrictive constraint.

For any $t \ge 0$, the random vector $X(t)$ will often be described in terms of the total number of customers $L(t)$ and the proportions of customers at the different nodes $\chi(t) = (\chi_i(t), 1 \le i \le n)$. More formally, define

$$L(t) = \sum_{j=1}^{n} X_j(t) = |X(t)| \quad \text{and} \quad \chi_i(t) = \frac{X_i(t)}{L(t)}, \qquad 1 \le i \le n, t \ge 0,$$

with the convention that $\chi(t) = e_1$ when $L(t) = 0$. Here, and more generally for any $x = (x_1, \dots, x_n) \in \mathbb{R}^n$, $|x|$ denotes the $\ell^1$ norm in $\mathbb{R}^n : |x| = \sum_1^n |x_i|$.

The vector $\chi(t)$ can be identified with a probability measure on $\{1, \dots, n\}$: namely, the empirical distribution of the positions of the $L(t)$ customers present in the network at time $t$. Denote by

$$\mathcal{P} = \left\{ \rho \in [0, +\infty[^n : \sum_{i=1}^{n} \rho_i = 1 \right\}$$

the state space of $\chi(t)$. The interior set of $\mathcal{P}$ is $\mathring{\mathcal{P}} = \{\rho \in ]0, +\infty[^n : \sum_1^n \rho_i = 1\}$.

As emphasized earlier, the deviation of $\chi(t)$ from $\pi$ will be of particular interest in the forthcoming analysis. It will be measured, depending on circumstances, by the $\ell^\infty$ distance $\|\chi(t) - \pi\|$,

$$\|x\| = \max_{1 \le i \le n} |x_i|, \qquad x = (x_1, \dots, x_n) \in \mathbb{R}^n,$$

or by the relative entropy $H(\chi(t), \pi)$ where $H(\cdot, \pi)$ is defined on the set $\mathcal{P}$ of probability measures on $\{1, \dots, n\}$ by

$$H(\rho, \pi) = \sum_{i=1}^{n} \rho_i \log \frac{\rho_i}{\pi_i} \in [0, +\infty[, \qquad \rho \in \mathcal{P}.$$

For $t \ge 0$, the quantity $H(\chi(t), \pi)$ will also be more simply denoted $H(t)$. The process $(H(t), t \ge 0)$ will spontaneously appear in the expression of the key martingale $J_\alpha$ introduced in the next section.

The different deviation times of $\chi(t)$ from $\pi$, or conversely, the time needed for $\chi(t)$ to reach a given neighborhood of $\pi$, will be of particular interest.

For $\varepsilon > 0$, $T_\varepsilon$ (resp. $T^\varepsilon$) denotes the first time when the $\ell^\infty$ distance between $\chi(t)$ and $\pi$ is smaller (resp. larger) than $\varepsilon$,

$$T_\varepsilon = \inf\{t \geq 0 : \|\chi(t) - \pi\| \leq \varepsilon\} \quad \text{and} \quad T^\varepsilon = \inf\{t \geq 0 : \|\chi(t) - \pi\| > \varepsilon\}.$$

Most results will be written in terms of these two stopping times, but it will be sometimes more convenient to work with the deviation time $T_H^\varepsilon$ from $\pi$ in terms of the relative entropy,

$$T_H^\varepsilon = \inf\{t \geq 0 : H(t) > \varepsilon\}.$$

All results on deviation times of $\chi(t)$ from $\pi$ defined in terms of the $\ell^\infty$ distance $\|\chi(t) - \pi\|$ can be translated into analogous estimates in terms of the relative entropy $H(t)$ thanks to the following classical result:

LEMMA 2.1. *There exist two $\pi$-depending positive constants $C_1$ and $C_2$ such that, for all $\rho \in \mathcal{P}$,*

$$C_1\|\rho - \pi\|^2 \leq H(\rho, \pi) \leq C_2\|\rho - \pi\|^2.$$

*In particular, for any $\varepsilon > 0$, $T_H^{C_1\varepsilon^2} \leq T^\varepsilon \leq T_H^{C_2\varepsilon^2}$.*

Another stopping time will play a central role, namely, the first time, denoted by $\mathcal{T}_0$, when the system has an empty node. Formally,

$$\mathcal{T}_0 = \inf\{t \geq 0 : \exists i \in \{1, \dots, n\}, X_i(t) = 0\}.$$

Indeed, the martingale property for the family of integrals presented in Section 3 will hold only up to time $\mathcal{T}_0$, that is, as long as the output rate at each node $i$ is exactly equal to $\mu_i$. In the same way, it will be easily shown that, for $t < \mathcal{T}_0$, $L(t)$ behaves exactly like the $M/M/1$ queue with input rate $\lambda$ and output rate $\mu$.

A last useful remark concerning these stopping times is that, when $\mathcal{T}_0$ is finite, $\|\chi(\mathcal{T}_0) - \pi\| \geq \min \pi_i (> 0)$. Together with Lemma 2.1, this immediately gives the following result:

LEMMA 2.2. *There exists $\varepsilon_0 > 0$ such that $T^\varepsilon \vee T_H^\varepsilon \leq \mathcal{T}_0$ holds for any $\varepsilon \leq \varepsilon_0$.*

**3. Martingale.** The results of this section are twofold: Theorem 3.1 gives the (almost) explicit expression of a local martingale $J_\alpha(\cdot \wedge \mathcal{T}_0)$, indexed by some positive parameter $\alpha$, and Proposition 3.2 derives the main estimate on deviation times $T_H^\varepsilon$ of $\chi(t)$ from $\pi$, that will be used in Sections 6 and 7. Concerning the construction of $J_\alpha$, the present section only aims at giving the main lines. The (numerous) technical details are postponed to the Appendix.

The approach for constructing the martingale $J_\alpha$ is similar to the approach used in [10] for the $M/M/\infty$ queue. The idea is to first exhibit a family of space–time harmonic functions $(h_v(t,x), v\in\mathbb{R}^n)$ for the generator $\Omega$ given by (1), and then to integrate $h_v(t,x)f(v)$ with respect to $v$ for some suitable function $f$, on some well-chosen, time-dependent domain. The last step is then to make a change of variables so that the new harmonic function is split into two factors, respectively, depending on time and space. The resulting local martingale is then adapted for an optional stopping use, leading to hitting-times estimations.

Some notation are required at this point. Denote by $(P_t, t\in\mathbb{R})$ the $Q$-generated Markov semi-group of linear operators in $\mathbb{R}^n: P_t = e^{tQ}$, extended to all real indices $t$ into a group. For $v\in\mathbb{R}^n$ and $t\in\mathbb{R}$, define

$$\phi(v,t) = (\phi_i(v,t), 1\le i\le n) = P_{-t}v.$$

Theorem 3.1 below requires the technical assumption that $Q$ is diagonalizable. Let $\theta$ be the trace of $-Q$, so that $\theta>0$, and let $\mathcal{S}\subset\mathbb{R}^{n-1}$ be the projection of $\mathring{\mathcal{P}}\subset\mathbb{R}^n$ on the $n-1$ first coordinates, that is,

$$\mathcal{S} = \left\{u = (u_1,\dots,u_{n-1})\in\mathbb{R}^{n-1}: \forall i = 1,\dots,n-1, u_i>0 \text{ and } \sum_{i=1}^{n-1} u_i < 1\right\}.$$

For any $u\in\mathcal{S}$, denote by $\tilde{u}\in\mathring{\mathcal{P}}$ the $n$th-dimensional vector which completes $u$ into a probability distribution, that is, $\tilde{u}_i = u_i$ for any $1\le i\le n-1$ and $\tilde{u}_n = 1-\sum_1^{n-1} u_i$.

The following proposition describes a family of space–time harmonic functions.

Proposition 3.1. *Let $v\in\mathbb{R}^n$ be fixed, and let $\varphi(v,\cdot)$ be any primitive of*

$$\sum_{i=1}^n\left(\mu_i\frac{\phi_i(v,\cdot)}{1+\phi_i(v,\cdot)} - \lambda_i\phi_i(v,\cdot)\right)$$

*on any open subset $V$ of $\{t\ge 0: 1+\phi_i(v,t)\ne 0$ for $i=1,\dots,n\}$. The function*

$$h_v(t,x) = e^{\varphi(v,t)}\prod_{i=1}^n (1+\phi_i(v,t))^{x_i}, \qquad t\in V, x\in\mathbb{N}^n,$$

*is space–time harmonic with respect to $\Omega$ in the domain $V\times\mathbb{N}^{*n}$.*

Proof. It must be shown that $\partial h_v(t,x)/\partial t + \Omega(h_v(t,\cdot))(x) = 0$ on the above domain. For $x\in\mathbb{N}^{*n}$ and $t\in V$, $h_v(t,x)\ne 0$, and one easily computes

$$\frac{1}{h_v(t,x)}\frac{\partial h_v}{\partial t}(t,x) = \frac{\partial\varphi}{\partial t}(v,t) + \sum_{i=1}^n x_i\frac{\partial\phi_i(v,t)/\partial t}{1+\phi_i(v,t)}$$

and

$$\frac{1}{h_v(t,x)}\Omega(h_v(t,\cdot))(x) = \sum_{i=1}^{n} \lambda_i \phi_i(v,t) - \sum_{i=1}^{n} \mu_i \frac{\phi_i(v,t)}{1+\phi_i(v,t)} + \sum_{1\le i\ne j\le n} x_i q_{ij} \frac{\phi_j(v,t)-\phi_i(v,t)}{1+\phi_i(v,t)}.$$

The last term in the right-hand side is equal to $\sum_{i=1}^{n} \frac{x_i}{1+\phi_i(v,t)}(Q\phi(v,t))_i$. By definition $\phi$ satisfies $\partial\phi(v,t)/\partial t = -Q\phi(v,t)$, and the result follows. □

REMARK 3.1. The product form of these space–time harmonic functions is quite similar to that of the harmonic functions introduced in [10] for the $M/M/\infty$ queue.

In addition, it is easily checked that, choosing $v = (u-1,\dots,u-1)$ for some $u \ne 0$, so that $v$ is some eigenvector of $P_t$, $t \in \mathbb{R}$, associated with eigenvalue 1, yields $h_v(t, X(t)) = u^{L(t)} e^{[\lambda(1-u)+\mu(1-1/u)]t}$, which is the martingale associated with an $M/M/1$ queue $L$ with arrival rate $\lambda$ and service rate $\mu$ (see, for example, [12]).

Starting from $h_v(t,x)$, two steps lead to $J_\alpha$: (i) integration of $h_v(t,x)$ over $v$ against some function $f(v)$ on a suitable time-dependent domain $\mathcal{D}(t)$; (ii) change of variables. These two steps are detailed and justified in the Appendix, yielding the following family of local martingales:

THEOREM 3.1. *There exist two positive, continuous, bounded functions $F$ and $G$ on $\mathring{\mathcal{P}}$ such that for any $\alpha > 0$, $u \mapsto F(\tilde{u})^{\alpha-1}$ is integrable on $\mathcal{S}$ and $(J_\alpha(t \wedge \mathcal{T}_0), t \ge 0)$ is a nonnegative local martingale where $J_\alpha(t)$ is defined for $\alpha > 0$ and $t \ge 0$ by*

$$J_\alpha(t) = e^{-\alpha\theta t} \int_{\mathcal{S}} \prod_{i=1}^{n} \left(\frac{\tilde{u}_i}{\pi_i}\right)^{X_i(t)} G(\tilde{u}) F(\tilde{u})^{\alpha-1}\, du,$$

*or equivalently,*

$$J_\alpha(t) = e^{-\alpha\theta t} \int_{\mathcal{S}} e^{L(t)(H(t)-H(\chi(t),\tilde{u}))} G(\tilde{u}) F(\tilde{u})^{\alpha-1}\, du. \tag{2}$$

*Moreover, $F$ satisfies*

$$\sup_{0<\alpha\le 1} \left(\alpha^n \int_{\mathcal{S}} F(\tilde{u})^{\alpha-1}\, du\right) < +\infty. \tag{3}$$

The advantage of $J_\alpha(t)$ [as compared to $h_v(t, X(t))$], is that the dependence in time is there split into two factors: $e^{-\alpha\theta t}$ is a direct function of time, and the integral is a function of the state of the system at time $t$, $X(t)$ or equivalently $(L(t), \chi(t))$.

The next proposition gives the fundamental estimate obtained through optional stopping and is used several times throughout the paper.

PROPOSITION 3.2. *For any $\delta$ such that $0 < \delta < \varepsilon_0$, where $\varepsilon_0$ is given by Lemma 2.2, there exists some constant $C_\delta$ such that*

$$\mathbb{E}_x(e^{-\alpha\theta T_H^\varepsilon}; L(T_H^\varepsilon) \geq \ell) \leq C_\delta \alpha^{-n} e^{|x|H(x/|x|,\pi)-(\varepsilon-\delta)\ell}$$

*holds for any initial state $x \in \mathbb{N}^n$ and any $\varepsilon \in ]\delta, \varepsilon_0[$, $\ell > 0$ and $\alpha \in ]0, 1]$.*

Proposition 3.2 is derived from the two following lemmas by choosing $T = T_H^\varepsilon$ (so that, by Lemma 2.2, $T \wedge \mathcal{T}_0 = T$ when $\varepsilon < \varepsilon_0$). Note that only Lemma 3.1 uses the fact that $J_\alpha$ is a local martingale, whereas Lemma 3.2 stems directly from the expression of $J_\alpha$ provided by (2).

LEMMA 3.1. *There exists some constant $C_3 > 0$ such that, for any $\alpha \in ]0, 1]$, any initial state $x \in \mathbb{N}^n$ and any stopping time $T$, the following inequality holds:*

$$\mathbb{E}_x[J_\alpha(T \wedge \mathcal{T}_0)] \leq C_3 \alpha^{-n} e^{|x|H(x/|x|,\pi)}.$$

PROOF. Fix $\alpha \in ]0, 1]$ and $x \in \mathbb{N}^n$. Since $J_\alpha(\cdot \wedge \mathcal{T}_0)$ is a nonnegative local martingale, it is a supermartingale, and so is $(J_\alpha(t \wedge T \wedge \mathcal{T}_0), t \geq 0)$ by Doob's optional stopping theorem. In particular, for any $t \geq 0$,

$$\mathbb{E}_x[J_\alpha(0)] \geq \mathbb{E}_x[J_\alpha(t \wedge T \wedge \mathcal{T}_0)]$$

and Fatou's lemma gives

$$\begin{aligned}\mathbb{E}_x[J_\alpha(0)] &\geq \liminf_{t\to+\infty} \mathbb{E}_x[J_\alpha(t \wedge T \wedge \mathcal{T}_0)] \\ &\geq \mathbb{E}_x\Big[\liminf_{t\to+\infty} J_\alpha(t \wedge T \wedge \mathcal{T}_0)\Big] = \mathbb{E}_x[J_\alpha(T \wedge \mathcal{T}_0)]\end{aligned}$$

[here $J_\alpha(T \wedge \mathcal{T}_0)$ makes sense a.s. when $T \wedge \mathcal{T}_0 = +\infty$ since any nonnegative supermartingale almost surely converges to some variable at infinity].

From the definition of $J_\alpha$ given by (2), using $e^{-y} \leq 1$ for $y \geq 0$, one gets

$$\mathbb{E}_x[J_\alpha(0)] \leq \sup_{\mathring{\mathcal{P}}}(G) e^{|x|H(x/|x|,\pi)} \int_{\mathcal{S}} F(\tilde{u})^{\alpha-1}\, du \leq C_3 e^{|x|H(x/|x|,\pi)} \alpha^{-n},$$

where $C_3 = \sup_{\mathring{\mathcal{P}}}(G) \sup_{\alpha\leq 1}(\alpha^n \int_{\mathcal{S}} F(\tilde{u})^{\alpha-1}\, du)$ is finite by (3), which proves the lemma. □

LEMMA 3.2. *For any positive $\delta$, there exists some positive constant $B_\delta$ such that the following implication holds for any $\alpha \in ]0,1]$, $\ell > 0$, $\varepsilon > \delta$ and $t \geq 0$:*

$$L(t) \geq \ell \text{ and } H(t) \geq \varepsilon \quad \Longrightarrow \quad J_\alpha(t) \geq B_\delta \cdot e^{-\alpha\theta t + (\varepsilon-\delta)\ell}.$$

PROOF. Fix $\varepsilon > \delta$, $\alpha \in ]0,1]$, $\ell > 0$ and $t \geq 0$. A lower bound on the integral part of (2) is obtained when $L(t) \geq \ell$ and $H(t) \geq \varepsilon$. For $v \in \mathcal{P}$, define the set $\mathcal{S}_\delta(v) \subset \mathcal{S}$ by

$$\mathcal{S}_\delta(v) = \{u \in \mathcal{S} : H(v, \tilde{u}) \leq \delta\}.$$

If $H(t) \geq \varepsilon$ and $L(t) \geq \ell$, then

$$\int_{\mathcal{S}} e^{L(t)(H(t)-H(\chi(t),\tilde{u}))} G(\tilde{u}) F(\tilde{u})^{\alpha-1}\, du \geq \beta e^{\ell(\varepsilon-\delta)} \int_{\mathcal{S}_\delta(\chi(t))} G(\tilde{u})\, du,$$

where $\beta = \min\{(\sup_{\overset{\circ}{\mathcal{P}}} F)^{-1}, 1\}$. Indeed, $\alpha$ being smaller than 1, $\beta$ is a lower bound for $F(\tilde{u})^{\alpha-1}$ on $\mathcal{S}$. Consider now the function $\Phi_\delta : \mathcal{P} \to \mathbb{R}^+$ defined by

$$\Phi_\delta(v) = \int_{\mathcal{S}_\delta(v)} G(\tilde{u})\, du.$$

Since $G$ is bounded, $\Phi_\delta$ is easily shown to be continuous (using, e.g., Lebesgue's theorem). Moreover, $\Phi_\delta(v) > 0$ for any $v \in \mathcal{P}$ [because $G > 0$ and the interior of $\mathcal{S}_\delta(v)$ is not empty], and since $\mathcal{P}$ is compact, $\inf_{\mathcal{P}} \Phi_\delta > 0$. Setting $B_\delta = \beta \inf_{\mathcal{P}} \Phi_\delta$ achieves the proof. □

**4. Two key representations.** The Markov process $(X(t), t \geq 0)$ with infinitesimal generator $\Omega$ defined by (1) can be seen as a particle system involving three types of transitions: births, deaths and migrations of particles from one site to another. The main purpose of this section is to show that $X$ can be decomposed into the difference of two pure birth and migration processes, up to some reflection term (Theorem 4.1). A simpler result (Proposition 4.1) tells that, as long as $X$ does not hit the axis, the process $L$ of the total number of particles just behaves as a random walk (or equivalently as an $M/M/1$ queue). Finally, a representation of process $X$ involving labelled particles is given in the case of null death rates.

Theorem 4.1, together with the latter representation, will be crucial for describing the unstable regime in Section 6, while Proposition 4.1 will be repeatedly used in the study of both the super and subcritical regimes.

The idea for decomposing $X$ is the following: when $\mu = 0$, the system consists of immortal particles generated at rate $\lambda$ and performing independent Markov trajectories. Introducing a death procedure, that is, some positive $\mu$, amounts to eliminating particles (at rate $\mu_i$ at site $i$) if possible. Up to some correction, due to the fact that no death can actually occur at an empty site, this is equivalent to subtracting some analogous process with birth rates $\mu_i$ $(1 \le i \le n)$, zero death rates and migration rate matrix $Q$.

This can be formalized by introducing an enlarged Markov process involving three types of particles. Define $(X, Y, Z)$ as a Markov process in $\mathbb{N}^{3n}$ with generator $\Gamma$ characterized by the following transitions and rates: for any $(x, y, z) \in \mathbb{N}^{3n}$ and $i, j \in \{1, \dots, n\}$ such that $i \ne j$,

$$(x,y,z) \longrightarrow \begin{cases} (x+e_i, y, z), & \text{at rate } \lambda_i, \\ (x-e_i, y+e_i, z), & \mu_i \mathbb{1}_{\{x_i \ge 1\}}, \\ (x, y, z+e_i), & \mu_i \mathbb{1}_{\{x_i = 0\}}, \\ (x-e_i+e_j, y, z), & q_{ij} x_i, \\ (x, y-e_i+e_j, z), & q_{ij} y_i, \\ (x, y, z-e_i+e_j), & q_{ij} z_i. \end{cases}$$

($X$ keeps track of the "real" particles, $Y$ of the killed ones and $Z$ of virtual particles generated at some site when no particle has been found to be killed.)

It is clear from these transitions and rates that, indexing generator $\Omega$ by its birth and death rate vectors: $\underline{\lambda} = (\lambda_i, 1 \le i \le n)$ and $\underline{\mu} = (\mu_i, 1 \le i \le n)$ $(\underline{\lambda}, \underline{\mu} \in [0, +\infty[^n)$ and denoting by $\underline{0}$ the null vector in $\mathbb{R}^n$:

(i) $X$ is a Markov process in $\mathbb{N}^n$ with generator $\Omega_{\underline{\lambda},\underline{\mu}}$;
(ii) $X + Y$ is also Markov in $\mathbb{N}^n$, with generator $\Omega_{\underline{\lambda},\underline{0}}$;
(iii) $Y + Z$ is Markov in $\mathbb{N}^n$ with generator $\Omega_{\underline{\mu},\underline{0}}$;
(iv) $|X + Y| - (|X(0) + Y(0)|)$ is some Poisson process with intensity $\lambda$;
(v) $|Y + Z| - (|Y(0) + Z(0)|)$ is some Poisson process with intensity $\mu$;
(vi) these two Poisson processes are independent.

Now from (i), any process $X$ with generator $\Omega$ can be considered as the first component of some Markov process with generator $\Gamma$ and with initial state $(X(0), \underline{0}, \underline{0})$.

The two next results are easily derived from this construction and from remarks (i) to (vi). In order to state the main theorem, it is convenient to index the process $X$ both by its initial state and by its birth and death parameters, writing $X^x_{\underline{\lambda},\underline{\mu}}$ for the process $X$ with initial state $x \in \mathbb{N}^n$, migration rate matrix $Q$ and birth (resp. death) parameters $\underline{\lambda} = (\lambda_i, 1 \le i \le n)$ [resp. $\underline{\mu} = (\mu_i, 1 \le i \le n)$].

THEOREM 4.1. *For any $x \in \mathbb{N}^n$ and $\underline{\lambda}, \underline{\mu} \in [0, +\infty[^n$, there exist versions of $X^x_{\underline{\lambda},\underline{\mu}}$, $X^x_{\underline{\lambda},\underline{0}}$ and $X^{\underline{0}}_{\underline{\mu},\underline{0}}$ such that*

$$X^x_{\underline{\lambda},\underline{\mu}} = X^x_{\underline{\lambda},\underline{0}} - X^{\underline{0}}_{\underline{\mu},\underline{0}} + Z,$$

*where $Z$ is an $\mathbb{N}^n$-valued process such that $|Z|$ is nondecreasing, initially zero, and increases only at times when some $X_i(t)$ is zero.*

PROOF. Write $X = X + Y - (Y + Z) + Z$, where $(X(0), Y(0), Z(0)) = (x, \underline{0}, \underline{0})$ and $(X, Y, Z)$ is Markov with generator $\Gamma$, so that $X$ is some version of $X^x_{\underline{\lambda},\underline{\mu}}$, and by (ii) and (iii), $X + Y$ is some version of $X^x_{\underline{\lambda},\underline{0}}$ and $Y + Z$ some version of $X^{\underline{0}}_{\underline{\mu},\underline{0}}$.

The theorem is proved, since $|Z|$ has the stated properties as can be seen on $\Gamma$. □

REMARK 4.1. The process $Z$ in Theorem 4.1 appears as a reflection term; it guarantees that $X$ stays nonnegative, compensating by adding some virtual particle for a jump of $X^{\underline{0}}_{\underline{\mu},0}$ that would get some $X_i$ to the value $-1$.

However, contrary to usual multidimensional Skorokhod reflection terms, here, due to the movements of particles, components $Z_i$'s are not necessarily nondecreasing in time; only their sum is. Also, $Z_i$ can increase at times when $X_i$ is not zero.

Theorem 4.1 and its proof, together with properties (iv), (v) and (vi), give the following proposition which constitutes one of the key ingredients for deriving the fluid limits in Sections 6 and 7.

PROPOSITION 4.1. *For all $t \le \mathcal{T}_0$, the following equality holds:*

$$L(t) = L(0) + \mathcal{N}_\lambda(t) - \mathcal{N}_\mu(t),$$

*where $\mathcal{N}_\lambda$ and $\mathcal{N}_\mu$ are independent Poisson processes with respective intensities $\lambda$ and $\mu$. Moreover, $L(t) \ge L(0) + \mathcal{N}_\lambda(t) - \mathcal{N}_\mu(t)$ holds for any $t \ge 0$.*

We conclude this section with a representation of process $X^x_{\underline{\lambda},\underline{0}}$ that will notably be used in Section 6, in conjunction with Theorem 4.1, for analyzing the unstable regime. $X^x_{\underline{\lambda},\underline{0}}$ is here obtained as a function of a Poisson process with intensity $\lambda$ and a sequence of Markov processes with infinitesimal generator $Q$ (representing the trajectories of the successively generated particles).

More precisely, $X^x_{\underline{\lambda},\underline{0}}$ admits the following representation:

$$X^x_{\underline{\lambda},\underline{0}}(t) = \left(\sum_{k\ge 1} \mathbb{1}_{\{\xi_k(t-\sigma_k)=i, \sigma_k \le t\}}, 1 \le i \le n\right), \qquad t \ge 0, \tag{4}$$

where:

- $\sigma_k = 0$ for $1 \le k \le |x|$;
- $\mathcal{N}_\lambda = (\sigma_k, k \ge |x| + 1)$ is a Poisson process with parameter $\lambda$;
- $\xi_k$, $k \ge 1$, are Markov jump processes in $\{1, \dots, n\}$ with generator $Q$ and initial distribution:
  - $\sum_{i=1}^n (\lambda_i/\lambda)\delta_i$ for $k \ge |x| + 1$;

– $\delta_i$ for $x_i$ arbitrarily chosen indices $k \in \{1, \dots, |x|\}$ $(1 \le i \le n)$;

• $\mathcal{N}_\lambda$ and the $\xi_k$, $k \ge 1$, are mutually independent.

$(\xi_k, 1 \le k \le |x|)$ hold for the trajectories of the initial particles and $(\xi_k, k \ge |x| + 1)$ for those of the successive newborn particles. $\mathcal{N}_\lambda$ holds for the global birth process. For $k \ge 1$, particle $k$ is in the system from time $\sigma_k$ ($\sigma_k = 0$ for $k \le |x|$).

Similarly as in the previous construction, a formal proof of (4) can be provided by constructing $X^x_{\underline{\lambda,0}}$ as function of a more complete process (that also contains $\mathcal{N}_\lambda$ and the $\xi_k$'s, $k \ge 1$), characterized through its infinitesimal generator and describing the list of current positions of particles present in the system ordered according to their birth rank.

**5. The space renormalized process.** The stability property of the system for $\lambda < \mu$ will be derived in Section 7 from a fluid scaling analysis, that is, from the study of the space–time renormalized process

$$\overline{X}^x(t) = \frac{X^x(|x|t)}{|x|}, \qquad t \ge 0,$$

as $|x|$ goes to infinity where $X^x$ is the Markov process $X$ initiated at $x$. It will be underlain by the $M/M/1$ behavior of the total occupancy process $L$ (only valid as long as no $X_i$ is zero, hence the intricacies of the analysis).

The particular behavior of $\overline{X}^x$ at $t = 0^+$ will result from the short-term behavior of the *only space* renormalized process $\widehat{X}$, defined as the the family of processes

$$\widehat{X}^x(t) = \frac{X^x(t)}{|x|}, \qquad t \ge 0, \text{ for } x \in \mathbb{N}^n \setminus \{0\}.$$

[The simpler notation $\widehat{X}(t)$, where $\widehat{X}(t) = X(t)/|X(0)|$, will also be used in situations where $|X(0)|$ is clearly nonzero.]

As highlighted in the Introduction, this scaling is natural and analogous to the Kelly scaling for the $M/M/\infty$ queue. This analogy appears in Proposition 5.3 below, that states convergence of $\widehat{X}^x$ as $|x| \to +\infty$ to some dynamical system having $\pi$ as its limiting point. In particular, for large $|x|$, $\widehat{X}^x$ reaches any neighborhood of $\pi$ in a quasi-deterministic finite time. And this will show (Sections 6 and 7) that asymptotically, $\overline{X}^x$ is instantaneously at $\pi$.

The results of this section are quite standard, essentially based on law of large numbers principles. The simple underlying idea is that as far as $\widehat{X}^x$ is only observed over a finite time window, since the number $|x|$ of initial particles goes to infinity while the numbers of births and deaths within the given window remain of the order of 1 (time is not rescaled here), the initial particles asymptotically dominate the system and mostly stay alive all along the time window thus behaving as $|x|$ independent Markov processes with generator $Q$.

For the same reasons, the process $\widehat{X}$ is not different, in the limit $|x| \to +\infty$, from the process $\chi = X/L$ of the spatial distribution of particles. The same convergence results hold for both processes; once proved for $\widehat{X}$, they easily extend to $\chi$.

Formalizing the above argument, the following coupling is intuitively clear. It compares the general model to the "closed" one (with no births nor deaths, but only initial particles). As in Section 4, generator $\Omega$ is indexed by its birth and death parameters $\underline{\lambda}$ and $\underline{\mu}$.

LEMMA 5.1. *For any $x \in \mathbb{N}^n$, there exists a coupling between the process $X^x$ with initial state $x$ and generator $\Omega_{\underline{\lambda},\underline{\mu}}$, and the process $U^x$ with initial state $x$ and generator $\Omega_{\underline{0},\underline{0}}$, such that for $t \geq 0$ and $i = 1, \dots, n$,*

$$U_i^x(t) - \mathcal{N}_\mu(t) \leq X_i^x(t) \leq U_i^x(t) + \mathcal{N}_\lambda(t),$$

*where $\mathcal{N}_\lambda$ and $\mathcal{N}_\mu$ are two Poisson processes with respective parameters $\lambda$ and $\mu$.*

*$X^x$, moreover, satisfies*

$$|x| - \mathcal{N}_\mu(t) \leq |X^x(t)| \leq |x| + \mathcal{N}_\lambda(t).$$

PROOF. The case $\mu = 0$ is a straightforward consequence of the representation (4) of $X$ from Section 4. Indeed if $\mu = 0$, (4) gives for $1 \leq i \leq n$,

$$U_i^x(t) \leq X_i^x(t) = \sum_{1 \leq k \leq |x|} \mathbb{1}_{\{\xi_k(t)=i\}} + \sum_{k \geq |x|+1} \mathbb{1}_{\{\xi_k(t-\sigma_k)=i, \sigma_k \leq t\}} \leq U_i^x(t) + \mathcal{N}_\lambda(t),$$

where $U^x(t)$ is constructed as $\sum_{k=1}^{|x|} \mathbb{1}_{\{\xi_k(t)=i\}}$ and $\mathcal{N}_\lambda = (\sigma_k, k \geq |x|+1)$.

Moreover, one gets $|U^x(t)| \leq |X^x(t)| \leq |U^x(t)| + \mathcal{N}_\lambda(t)$ by summing up over $i$ the previous first inequalities. The lemma is proved in this case since $|U^x(t)| = |x|$ for any $t \geq 0$.

The general case is then derived using the first part of Section 4. Indeed, consider $X^x$ as the first component of a random process $(X^x, Y, Z)$ in $\mathbb{N}^{3n}$ such that $Y(0) = Z(0) = 0$, $X^x + Y$ is some process with generator $\Omega_{\underline{\lambda},\underline{0}}$ and $|Y + Z|$ is some Poisson process $\mathcal{N}_\mu$ with intensity $\mu$. The first part of the proof then applies to $X^x + Y$ and gives, for $t \geq 0$,

$$U^x(t) \leq X^x(t) + Y(t) \leq U^x(t) + \mathcal{N}_\lambda(t)$$

componentwise, as well as $|x| \leq |X^x(t) + Y(t)| \leq |x| + \mathcal{N}_\lambda(t)$.

The lemma follows by noticing that

$$0 \leq Y_i(t) \leq |Y(t)| \leq |Y(t) + Z(t)| = \mathcal{N}_\mu(t), \qquad 1 \leq i \leq n. \qquad \square$$

The two main results of this section concern the hitting time of some neighborhood of $\pi$ by the space renormalized process $\widehat{X}$. Namely, for any positive $\delta$,

$$\widehat{T}_\delta = \inf\{t \geq 0 : \|\widehat{X}(t) - \pi\| \leq \delta\}.$$

Recall that the analogous time with $\chi$ in place of $\widehat{X}$ is denoted by $T_\delta$.

PROPOSITION 5.1. *For any positive $\delta$, there exists some deterministic time $t_\delta \geq 0$ such that*

$$\lim_{|x|\to+\infty} \mathbb{P}_x(\widehat{T}_\delta > t_\delta) = 0.$$

*The same result holds for the stopping time $T_\delta$.*

PROOF. We refer to the proof of Proposition 5.2 below. Proposition 5.1 is obtained in the same way, just changing $\delta_N$, $s_N$ and $t_N$ into $\delta$, $s = -1/\eta \log(\delta/2B)$ and $t_\delta = -1/\eta \log(\delta/4B)$. $\square$

The following more accurate result will be required for analyzing the subcritical case $\lambda < \mu$.

PROPOSITION 5.2. *There exist two positive constants $A$ and $\eta$ such that, for any sequence of positive numbers $(\delta_N, N \geq 1)$ satisfying*

$$\lim_{N\to+\infty} \delta_N = 0 \quad \textit{and} \quad \lim_{N\to+\infty} \delta_N \sqrt{N} = +\infty,$$

*then*

$$\lim_{N\to+\infty} \Big[\max_{x\in\mathbb{N}^n : |x|=N} \mathbb{P}_x(\widehat{T}_{\delta_N} > t_N)\Big] = 0 \qquad \textit{where } t_N = -\frac{1}{\eta}\log\frac{\delta_N}{A}.$$

*The same result holds for the stopping time $T_{\delta_N}$.*

PROOF. First consider a closed system, that is, assume $\lambda = \mu = 0$; the general case will then be deduced from Lemma 5.1. As in Lemma 5.1, let $U^x$ be the closed process with initial state $x \in \mathbb{N}^n$ where $|x| = N$. In this case (4) becomes

$$U_i^x(t) = \sum_{k=1}^N \mathbb{1}_{\{\xi_k(t)=i\}}, \qquad 1 \leq i \leq n, t \geq 0,$$

where $\xi_k$, $1 \leq k \leq N$, are independent Markov processes with the same generator $Q$ and different initial conditions: for any $i$, $1 \leq i \leq n$, $\xi_k(0) = i$ for $x_i$ of the $N$ indices $k = 1, \dots, N$.

As introduced in Section 3, let $(P_t, t \geq 0)$ denote the transition semi-group associated to $Q$. The exponentially fast convergence of any irreducible finite state space Markov semi-group to its stationary distribution, tells existence of $B > 0$ and $\eta > 0$ such that

$$\max_{1\leq i,j\leq n} |P_t(j,i) - \pi_i| \leq Be^{-\eta t}, \qquad t \geq 0.$$

In particular, $\max_{1\leq i,j\leq n}|P_{s_N}(j,i) - \pi_i| \leq \delta_N/2$ for $s_N = -\frac{1}{\eta}\log\frac{\delta_N}{2B}$.

The outline of the proof for the closed case is the following: at time $s_N$, all trajectories $\xi_k$, $1 \leq k \leq N$, are very close to $\pi$ in distribution (by the order of $\delta_N$). Since $\widehat{U}^x(t)$ represents the empirical distribution of the $N$ particles at time $t$, the law of large numbers shows that for large $N$, $\widehat{U}^x(s_N)$ is also close to $\pi$ (by the same order), because $\delta_N$ tends to 0 not too fast

Precisely, for any $N \geq 1$ and $x \in \mathbb{N}^n$ such that $|x| = N$,

$$\|\mathbb{E}(\widehat{U}^x(s_N)) - \pi\| = \left\|\mathbb{E}\left(\frac{U^x(s_N)}{N}\right) - \pi\right\| = \left\|\frac{1}{N}\sum_{k=1}^{N}(\mathbb{P}(\xi_k(s_N) = \cdot) - \pi)\right\| \leq \frac{\delta_N}{2}.$$

Thus, for any $N \geq 1$, using Chebyshev's inequality for the last step,

$$\begin{aligned}\mathbb{P}(\|\widehat{U}^x(s_N) - \pi\| > \delta_N) &\leq \mathbb{P}\left(\|\widehat{U}^x(s_N) - \mathbb{E}(\widehat{U}^x(s_N))\| > \frac{\delta_N}{2}\right)\\ &\leq \sum_{i=1}^{n}\mathbb{P}\left(|\widehat{U}_i^x(s_N) - \mathbb{E}(\widehat{U}_i^x(s_N))| > \frac{\delta_N}{2}\right)\\ &\leq \sum_{i=1}^{n}\frac{\mathrm{Var}(U_i^x(s_N))}{\delta_N^2 N^2/4}.\end{aligned}$$

Independence of the processes $(\xi_k, 1 \leq k \leq N)$ yields

$$\mathrm{Var}(U_i^x(s_N)) = \sum_{k=1}^{N}\mathrm{Var}(\mathbb{1}_{\{\xi_k(s_N)=i\}}) \leq \frac{N}{4}$$

(bounding the variance of any Bernoulli random variable by 1/4). Finally,

$$\max_{x\in\mathbb{N}^n\,:\,|x|=N}\mathbb{P}(\|\widehat{U}^x(s_N) - \pi\| > \delta_N) \leq \frac{n}{\delta_N^2 N}. \tag{5}$$

Now consider the process $X^x$ associated to any family $(\lambda_i, \mu_i, 1 \leq i \leq n)$ of parameters and any initial state $x$ such that $|x| = N$. Still denote by $U^x$ the associated closed process with the same initial state $x$.

Define $t_N = -\frac{1}{\eta}\log\frac{\delta_N}{4B}$. The first part of Lemma 5.1 implies that, for any $N \geq 1$,

$$\|\widehat{X}^x(t_N) - \pi\| \leq \|\widehat{U}^x(t_N) - \pi\| + \left\|\frac{1}{N}(\mathcal{N}_\lambda(t_N) + \mathcal{N}_\mu(t_N))\right\|,$$

so that

$$\mathbb{P}_x(\widehat{T}_{\delta_N} > t_N) \le \mathbb{P}\left(\|\widehat{U}^x(t_N) - \pi\| > \frac{\delta_N}{2}\right) + \mathbb{P}\left(\|\mathcal{N}_\lambda(t_N) + \mathcal{N}_\mu(t_N)\| > \frac{N\delta_N}{2}\right).$$

By (5) the first term tends to zero uniformly in $x$ as $N$ goes to infinity, since $t_N$ is associated to $\delta_N/2$ in the same way as $s_N$ was to $\delta_N$. The second one is also easily shown to converge to zero, using Chebyshev's inequality for the Poisson variable $\mathcal{N}_\lambda(t_N) + \mathcal{N}_\mu(t_N)$, together with the relation $\delta_N\sqrt{N} \gg 1$ that implies $N\delta_N \gg \sqrt{N} \gg 1/\delta_N \gg t_N$.

The first part of the proposition is thus proved with $A = 4B$.

Using the last assertion of Lemma 5.1, it is not difficult to show that the same result holds for $T_{\delta_N}$. □

We finally just mention for the sake of completeness (it will not be used in the sequel) the following result that describes the asymptotic dynamics, as $|x| \to +\infty$, of the empirical distribution of the particles; it evolves as the distribution, as function of time, of a Markov process with generator $Q$.

Not surprisingly, this can be proved using the same standard arguments as for studying the $M/M/\infty$ queue under the Kelly scaling (see [12]).

PROPOSITION 5.3. *Consider the processes* $(\widehat{X}^{x_N}(t), t \ge 0)$ *associated with some sequence* $(x_N, N \ge 1)$ *of initial states satisfying* $\lim_{N\to+\infty} \frac{x_N}{N} = \rho$, *for some* $\rho \in \mathcal{P}$.

*For any* $T > 0$, *as* $N \to +\infty$, $(\widehat{X}^{x_N}(t), t \ge 0)$ *converges in distribution with respect to the uniform norm topology on* $[0, T]$, *to the deterministic trajectory*

$$\rho(t) = \rho P_t.$$

*In other words, for any positive* $\delta$,

$$\lim_{N\to+\infty} \mathbb{P}\left(\sup_{0\le t\le T} \|\widehat{X}^{x_N}(t) - \rho P_t\| > \delta\right) = 0.$$

*The same convergence holds for the corresponding processes* $(\chi^{x_N}(t), t \ge 0)$, $N \ge 1$.

**6. The supercritical regime.** This section deals with the supercritical regime $\lambda > \mu$. As the next proposition shows, the instability of the system is straightforward in this case. Theorem 6.1 establishes an almost sure result describing the long-term behavior, and Theorem 6.2 presents a surprising phenomenon.

PROPOSITION 6.1. *When* $\lambda > \mu$ *the process* $X$ *is not ergodic.*

PROOF. Just remark, using Proposition 4.1, that if $x \in \mathbb{N}^n$ is the initial state,

$$L(t) \geq |x| + \mathcal{N}_\lambda(t) - \mathcal{N}_\mu(t).$$

Hence for any initial state, $L(t)$ almost surely goes to $+\infty$ as $t$ tends to $+\infty$. $\square$

The following theorem gives an almost sure description of the divergence of $X(t)$ for $t$ large. Among other arguments, the proof makes use for the first time of the martingale estimate provided by Proposition 3.2, and involves the representations of $X$ given in Section 4.

THEOREM 6.1. *Assume $\lambda > \mu$. Then, for any initial state $x \in \mathbb{N}^n$, the following convergence holds almost surely:*

$$\lim_{t\to+\infty} \frac{X^x(t)}{t} = (\lambda - \mu)\pi.$$

REMARK 6.1. This theorem has a double meaning: it tells almost sure convergence both of $\chi(t)$ to $\pi$ and of $L(t)/t$ to $\lambda - \mu$ as $t \to +\infty$.

PROOF OF THEOREM 6.1. Assume the theorem is true when $\mu = 0$. Then, using the notation of Theorem 4.1, $t^{-1}(X^x_{\underline{\lambda},\underline{0}}(t) - X^{\underline{0}}_{\underline{\mu},\underline{0}}(t))$ converges a.s. to $(\lambda-\mu)\pi$ and the componentwise inequality $X^x_{\underline{\lambda},\underline{\mu}} \geq X^x_{\underline{\lambda},\underline{0}} - X^{\underline{0}}_{\underline{\mu},\underline{0}}$ derived from Theorem 4.1, implies that each $X^x_i(t)$ tends to infinity almost surely as $t$ goes to infinity.

As a consequence, since $|Z(t)|$ can increase only when some $X_i(t)$ is zero, then, with probability 1, $\lim_{t\to+\infty} |Z(t)|$ is finite and $\lim_{t\to+\infty} Z(t)/t = 0$, so that

$$\lim_{t\to+\infty} \frac{X^x(t)}{t} = \lim_{t\to+\infty} \frac{X^x_{\underline{\lambda},\underline{0}}(t) - X^{\underline{0}}_{\underline{\mu},\underline{0}}(t)}{t} = (\lambda - \mu)\pi$$

holds almost surely, which is the stated result.

The theorem must now be proved in the case where $\mu = 0$. In this case with no deaths, using representation (4), the process $X^x$ splits into two (independent) processes: $X^x = U^x + X^{\underline{0}}$, where $U^x$ is associated to a "closed" system with $|x|$ particles moving independently, and $X^{\underline{0}}$ has no initial particles, birth rates $\underline{\lambda}$ and null death rates. Then $t^{-1}U^x(t)$ obviously tends to zero almost surely as $t$ tends to infinity, and all that is left to show is that $t^{-1}X^{\underline{0}}(t)$ converges almost surely to $\lambda\pi$.

So dropping for simplicity the superscript $\underline{0}$, consider the process $X$ with initial state $\underline{0}$, birth rates $\underline{\lambda}$ and null death rates. Equation (4) here becomes for $t \geq 0$,

$$X_i(t) = \sum_{k=1}^{\mathcal{N}_\lambda(t)} \mathbb{1}_{\{\xi_k(t-\sigma_k)=i\}}, \qquad 1 \leq i \leq n,$$

where $(\xi_k, k \geq 1)$ have initial distribution $\sum_{i=1}^n (\lambda_i/\lambda)\delta_i$.

It will first be shown that the analysis can be reduced to the case of *stationary* trajectories (i.e., the case when $\lambda_i/\lambda = \pi_i$ for $1 \leq i \leq n$) by using a coupling argument.

Indeed, associate with each $\xi_k$ a stationary process $\xi'_k$ with the same generator, such that $((\xi_k, \xi'_k), k \geq 1)$ is a sequence of independent processes in $\{1, \ldots, n\}^2$, and, for $k \geq 1$, $\xi_k$, $\xi'_k$ are coupled in the classical following way: $\xi_k$ and $\xi'_k$ are independent until the first time $T_k$ when they meet, and after that stay equal for ever. Recall that the "coupling times" $T_k$, $k \geq 1$, are integrable. Moreover, assume the $(\xi_k, \xi'_k), k \geq 1$, independent from $\mathcal{N}_\lambda$.

Define the process $(X'(t), t \geq 0)$ on $\mathbb{N}^n$ analogously to $X$, with the same $\mathcal{N}_\lambda$, but with $\xi'_k$ in place of $\xi_k$ $(k \geq 1)$. Then, for each $i \in \{1, \ldots, n\}$,

$$|X_i(t) - X'_i(t)| = \left| \sum_{k=1}^{\mathcal{N}_\lambda(t)} (\mathbb{1}_{\{\xi_k(t-\sigma_k)=i\}} - \mathbb{1}_{\{\xi'_k(t-\sigma_k)=i\}}) \right| \leq \sum_{k=1}^{\mathcal{N}_\lambda(t)} \mathbb{1}_{\{T_k > t-\sigma_k\}}.$$

Denoting $A(t)$ the last term, $A(t)$ is exactly the number of customers at time $t$ in an $M/G/\infty$ queue with no customer at time 0, arrival process $\mathcal{N}_\lambda$, and services given by the i.i.d. integrable variables $T_k$, $k \geq 1$.

It is easily proved that $A(t)/t$ converges almost surely to zero as $t$ tends to infinity. It is then enough to prove a.s. convergence of process $X'$ to $\lambda\pi$, and so we assume from now on that $(\xi_k, k \geq 1)$ are stationary.

Since $L(t)/t = \mathcal{N}_\lambda(t)/t$ converges a.s. to $\lambda$ as $t$ tends to infinity, the problem is equivalent to proving that $\chi(t)$ converges almost surely to $\pi$, that is, by Lemma 2.1 that

$$\forall \varepsilon > 0 \qquad \mathbb{P}(\exists T < +\infty : \forall t \geq T, H(t) \leq \varepsilon) = 1.$$

This will be done using Borel–Cantelli lemma and showing that

$$\forall \varepsilon > 0 \qquad \sum_{k=1}^{+\infty} \mathbb{P}(\exists t \in [\sigma_k, \sigma_{k+1}[ : H(t) > \varepsilon) < +\infty.$$

Writing, for any fixed $\varepsilon$,

$$\begin{aligned} &\mathbb{P}(\exists t \in [\sigma_k, \sigma_{k+1}[ : H(t) > \varepsilon) \\ &\quad \leq \mathbb{P}\left(H(\sigma_k) > \frac{\varepsilon}{2}\right) + \mathbb{P}\left(H(\sigma_k) \leq \frac{\varepsilon}{2} \text{ and } \exists t \in ]\sigma_k, \sigma_{k+1}[ : H(t) > \varepsilon\right), \end{aligned} \tag{6}$$

we will show that both series associated with both terms in the right-hand side converge for $\varepsilon$ sufficiently small [which is enough by monotonicity of the left-hand side of (6)].

Let us begin with the first term. Note that for $k \geq 1$, $\chi(\sigma_k) = X(\sigma_k)/k$. Then due to Lemma 2.1, it is enough to show that, for small $\varepsilon$ and any $i \in \{1, \dots, n\}$,

$$\sum_{k=1}^{+\infty} \mathbb{P}\left(\left|\frac{X_i(\sigma_k)}{k} - \pi_i\right| > \varepsilon\right) < +\infty. \tag{7}$$

This is obtained by using Chernoff's inequality, that we recall in Lemma 6.1.

Lemma 6.1 (Chernoff's inequality). *Let $Z_h, 1 \leq h \leq k$, be $k$ independent random variables such that $|Z_h| \leq 1$ and $\mathbb{E}(Z_h) = 0$ for $1 \leq h \leq k$.*

*The following bound holds for any $\eta \in [0, 2\sigma]$ where $\sigma^2 = \operatorname{Var}(\sum_{h=1}^k Z_h)$:*

$$\mathbb{P}\left(\left|\sum_{h=1}^{k} Z_h\right| \geq \eta\sigma\right) \leq 2e^{-\eta^2/4}.$$

Write

$$\frac{X_i(\sigma_k)}{k} - \pi_i = \frac{1}{k}\sum_{h=1}^{k} Z_{k,h}^{(i)} \qquad \text{with } Z_{k,h}^{(i)} = \mathbb{1}_{\{\xi_h(\sigma_k - \sigma_h) = i\}} - \pi_i, 1 \leq i \leq n.$$

Since $\xi_h, h \geq 1$, are stationary, then for each fixed $i \in \{1, \dots, n\}$ and $k \geq 1$, the $k$ variables $Z_{k,h}^{(i)}$, $1 \leq h \leq k$, are i.i.d. centered random variables, bounded by 1 in modulus. (Notice that independence is only true in this stationary case.)

We can thus apply Chernoff's inequality, which gives, for each fixed $k$ and $i$,

$$\mathbb{P}\left(\left|\frac{X_i(\sigma_k)}{k} - \pi_i\right| > \varepsilon\right) = \mathbb{P}\left(\left|\sum_{h=1}^{k} Z_{k,h}^{(i)}\right| > k\varepsilon\right) \leq 2e^{-\varepsilon^2 k/(4v_i)},$$

if $\varepsilon \leq 2v_i$ where $v_i = \pi_i(1 - \pi_i)$ is the common variance of the variables $Z_{k,h}^{(i)}$. Property (7) is then proved (for small $\varepsilon$, hence for any $\varepsilon$ by monotonicity).

Now it must be shown that the second term in the right-hand side of (6) is summable as well for $\varepsilon$ small enough. Here, the stationarity of the movements will play no special role.

By definition of $\sigma_k$, $\chi(t) = X(t)/k$ for any $t \in [\sigma_k, \sigma_{k+1}[$. Moreover, $\sigma_k$ is a stopping time for the Markov process $(X(t), t \geq 0)$, because it is the first time when $L(t) = k$. Hence the strong Markov property yields

$$\mathbb{P}_{\underline{0}}\left(H(\sigma_k) \leq \frac{\varepsilon}{2} \text{ and } \exists t \in ]\sigma_k, \sigma_{k+1}[ : H(t) > \varepsilon\right) \leq \max_{\substack{x \in \mathbb{N}^n : |x| = k \text{ and} \\ H(x/|x|, \pi) \leq \varepsilon/2}} \mathbb{P}_x(T_H^\varepsilon < \sigma_1).$$

Clearly, the last event only depends on $\sigma_1$ and on the movements of the $|x|$ initial particles, so that by independence of these variable and processes, one obtains, for any $x \in \mathbb{N}^n$,

$$\mathbb{P}_x(T_H^\varepsilon < \sigma_1) = \mathbb{E}_x(e^{-\lambda \widetilde{T}_H^\varepsilon}) \leq \mathbb{E}_x(e^{-(\lambda \wedge \theta)\widetilde{T}_H^\varepsilon}),$$

where $\widetilde{T}_H^\varepsilon$ is the first time the entropy associated to the initial particles is larger than $\varepsilon$. Then using Proposition 3.2 in the case of a closed system with $\delta = \varepsilon/4$, $\alpha = (\lambda/\theta) \wedge 1$ and $\ell = k$ gives

$$\mathbb{P}_x(T_H^\varepsilon < \sigma_1) \leq C_{\varepsilon/4}[(\lambda/\theta) \wedge 1]^{-n} e^{-\varepsilon k/4}$$

for any $x \in \mathbb{N}^n$ such that $|x| = k$ and $H(x/|x|, \pi) \leq \varepsilon/2$. The second term in (6) is thus summable over $k$ for $\varepsilon$ small enough. □

Along the preceding proof, we used $\sigma_1$, in the particular case $\mu = 0$, as an asymptotic lower bound (as the initial state grows to infinity) for the exit time of $\chi(t)$ from some neighborhood of $\pi$. This is a very crude underestimation, as the following result shows that this exit time is actually infinite with high probability.

THEOREM 6.2. *Assume $\lambda > \mu$, and fix $\delta$ and $\varepsilon$ such that $0 < \delta < \varepsilon < \varepsilon_0$ where $\varepsilon_0$ is given by Lemma 2.2. Consider a sequence $(x_N, N \geq 1)$ with $\lim_{N\to+\infty} |x_N|/N = 1$ and $H(x_N/|x_N|, \pi) \leq \delta$. Then*

$$\lim_{N\to+\infty} \mathbb{P}_{x_N}(T_H^\varepsilon = +\infty) = 1.$$

PROOF. By definition of $T_H^\varepsilon$, $\mathbb{P}_{x_N}(T_H^\varepsilon < +\infty) = \mathbb{P}_{x_N}(\exists t \geq 0 : H(t) \geq \varepsilon)$, and so we need to study the behavior of $H(t)$ for all time $t \geq 0$. The idea of the proof is twofold: first, the estimate given by Proposition 3.2 is precise enough to show that $T_H^\varepsilon$ is much larger than $N$, say $T_H^\varepsilon \geq N^2$. After this time, the initial particles are negligible, and Theorem 6.1 then gives a control on the rest of the trajectory by reducing the problem to the case where the system starts empty. So we use the following decomposition:

$$\mathbb{P}_{x_N}(T_H^\varepsilon < +\infty) \leq \mathbb{P}_{x_N}(T_H^\varepsilon \leq N^2) + \mathbb{P}_{x_N}(\exists t \geq N^2 : H(t) \geq \varepsilon).$$

For the first term, Markov's inequality gives

$$\mathbb{P}_{x_N}(T_H^\varepsilon \leq N^2) \leq e\mathbb{E}_{x_N}(e^{-T_H^\varepsilon/N^2}). \tag{8}$$

Let $\delta' < \varepsilon - \delta$. By choice of $\varepsilon$ and $\delta$, and since $H(x_N/|x_N|, \pi) \le \delta$, Proposition 3.2 shows that there exists a constant $C_{\delta'}$ such that by choosing $\alpha = 1/(\theta N^2)$, for any $N$ large enough and any $\ell_N$,

$$\mathbb{E}_{x_N}(e^{-T_H^\varepsilon/N^2}; L(T_H^\varepsilon) \ge \ell_N) \le C_{\delta'} e^{\delta|x_N| + 2n \log N - (\varepsilon - \delta')\ell_N}.$$

The choice of $\ell_N$ requires some care; as $N$ grows, it must be both of order $|x_N|$ and smaller than $L(T_H^\varepsilon)$ with high probability. Since $|x_N| \sim N$, write $|x_N| = N + u_N$ with $u_N = o(N)$, and choose $\ell_N = N - \sqrt{N v_N}$ with $v_N = |u_N| \vee 1$. With this choice, $\ell_N \sim N$ and $\ell_N - |x_N| \to -\infty$. The first relation implies, since $\varepsilon - \delta' - \delta > 0$,

$$\lim_{N \to +\infty} e^{\delta|x_N| + 2n \log N - (\varepsilon - \delta')\ell_N} = 0.$$

Moreover, since $T_H^\varepsilon \le \mathcal{T}_0$ because $\varepsilon < \varepsilon_0$, Proposition 4.1 implies that $L(T_H^\varepsilon) = L(0) + \mathcal{N}_\lambda(T_H^\varepsilon) - \mathcal{N}_\mu(T_H^\varepsilon)$, hence

$$\begin{aligned}\mathbb{P}_{x_N}(L(T_H^\varepsilon) \le \ell_N) &= \mathbb{P}_{x_N}(|x_N| + \mathcal{N}_\lambda(T_H^\varepsilon) - \mathcal{N}_\mu(T_H^\varepsilon) \le \ell_N) \\ &\le \mathbb{P}\Big(\inf_{t \ge 0}(\mathcal{N}_\lambda(t) - \mathcal{N}_\mu(t)) \le \ell_N - |x_N|\Big),\end{aligned}$$

where the last bound vanishes because $\lambda > \mu$, and so $\inf_{t \ge 0}(\mathcal{N}_\lambda(t) - \mathcal{N}_\mu(t))$ is finite with probability one whereas $\ell_N - |x_N|$ goes to $-\infty$. It results that $\mathbb{P}_{x_N}(T_H^\varepsilon \le N^2)$ goes to 0 thanks to (8) and to the following inequality:

$$\mathbb{E}_{x_N}(e^{-T_H^\varepsilon/N^2}) \le \mathbb{E}_{x_N}(e^{-T_H^\varepsilon/N^2}; L(T_H^\varepsilon) \ge \ell_N) + \mathbb{P}_{x_N}(L(T_H^\varepsilon) \le \ell_N)$$

and it has been shown that each term goes to 0.

All that is left to prove now is that $\lim_{N \to +\infty} \mathbb{P}_{x_N}(\exists t \ge N^2 : H(t) \ge \varepsilon) = 0$, or, by Lemma 2.1, that $\mathbb{P}_{x_N}(\exists t \ge N^2 : \|\chi(t) - \pi\| \ge \varepsilon)$ vanishes. After time $N^2$, the initial particles are negligible since a number of new particles of the order of $N^2$ have arrived. So the behavior of the system will be similar to that of a system starting empty, to which we can apply Theorem 6.1 (since in this case the initial state is fixed).

To formalize this argument, a coupling between the processes $X^x$ and $X^0$, for any $x \in \mathbb{N}^n$, is required.

LEMMA 6.2. *For any $x, y \in \mathbb{N}^n$ with $x \ge y$ componentwise, it is possible to couple the two processes $X^x$ and $X^y$ in such a way that for any $t \ge 0$, $L^x(t) - L^y(t) \le |x| - |y|$ and the inequality $X^x(t) \ge X^y(t)$ holds componentwise.*

The proof of this lemma is postponed at the end of the current proof. Let $X^0$ be the process starting empty coupled with $X^{x_N}$, and let $L^0 =$ $|X^0|$, $L^{x_N} = |X^{x_N}|$, $\chi^0 = X^0/L^0$ and $\chi^{x_N} = X^{x_N}/L^{x_N}$ be the corresponding quantities. The triangular inequality gives

$$\begin{aligned}(9)\qquad &\mathbb{P}(\exists t \ge N^2 : \|\chi^{x_N}(t) - \pi\| \ge \varepsilon)\\ &\quad \le \mathbb{P}(\exists t \ge N^2 : \|\chi^{x_N}(t) - \chi^0(t)\| \ge \varepsilon/2)\\ &\qquad + \mathbb{P}(\exists t \ge N^2 : \|\chi^0(t) - \pi\| \ge \varepsilon/2).\end{aligned}$$

Theorem 6.1 states that $\chi^0(t)$ converges to $\pi$ almost surely, which shows that the last term goes to 0. For the first term, write for each $i = 1, \dots, n$,

$$\chi_i^{x_N}(t) - \chi_i^0(t) = \frac{X_i^{x_N}(t)}{L^{x_N}(t)} - \frac{X_i^0(t)}{L^0(t)} = \frac{(X_i^{x_N}(t) - X_i^0(t))L^0(t) - X_i^0(t)\Delta^{x_N}(t)}{L^0(t)(\Delta^{x_N}(t) + L^0(t))},$$

where $\Delta^{x_N}(t) = L^{x_N}(t) - L^0(t)$. Lemma 6.2 implies that $|X_i^{x_N}(t) - X_i^0(t)| \le \Delta^{x_N}(t) \le |x_N|$, hence since the function $z \mapsto z/(z+a)$ is decreasing for any $a \ge 0$,

$$|\chi_i^{x_N}(t) - \chi_i^0(t)| \le \frac{2\Delta^{x_N}(t)}{\Delta^{x_N}(t) + L^0(t)} \le \frac{2|x_N|}{|x_N| + L^0(t)} = \frac{2}{1 + L^0(t)/|x_N|}.$$

This yields in turn, using $t \ge N^2$ for the second inequality,

$$\begin{aligned}&\mathbb{P}(\exists t \ge N^2 : \|\chi^{x_N}(t) - \chi^0(t)\| \ge \varepsilon/2)\\ &\quad \le \mathbb{P}\left(\exists t \ge N^2 : \frac{2}{1 + L^0(t)/|x_N|} \ge \varepsilon/2\right)\\ &\quad \le \mathbb{P}\left(\inf_{t \ge N^2}(1 + L^0(t)/t \cdot N^2/|x_N|) \le 4/\varepsilon\right).\end{aligned}$$

Theorem 6.1 shows that $L^0(t)/t \to \lambda - \mu$ almost surely as $t \to +\infty$, and $N^2/|x_N|$ goes to infinity as $N$ goes to infinity by choice of $x_N$. Hence almost surely,

$$\lim_{N \to +\infty} \inf_{t \ge N^2}(1 + L^0(t)/t \cdot N^2/|x_N|) = +\infty$$

and the theorem is proved. □

We now fill in the gap in this proof by proving Lemma 6.2.

Proof of Lemma 6.2. Process $X$ admits the following representation as the solution of a system of integral equations:

$$\begin{aligned}X_i(t) = X_i(0) + \mathcal{N}_{\lambda_i}(t) - \int_0^t \mathbb{1}_{\{X_i(s^-)\ge 1\}}\, d\mathcal{N}_{\mu_i}(s)\\ + \sum_{j \ne i}\int_0^t \sum_{k=1}^{X_j(s^-)} d\mathcal{N}_{q_{ji}}^k(s) - \sum_{j\ne i}\int_0^t \sum_{k=1}^{X_i(s^-)} d\mathcal{N}_{q_{ij}}^k(s), \qquad 1 \le i \le n,\end{aligned}$$

where $\mathcal{N}_{\lambda_i}$ and $\mathcal{N}_{\mu_i}$, for $i = 1, \dots, n$, are Poisson processes with respective parameters $\lambda_i$ and $\mu_i$, and for $(i,j) \in \{1,\dots,n\}^2$, $i \neq j$, $(\mathcal{N}^k_{q_{ij}}, k \geq 1)$ is a sequence of Poisson processes with parameter $q_{ij}$, all these processes being independent.

Now using the same Poisson processes for $X^x$ and $X^y$, it is easy to check that the inequalities $X_i^x(t) \geq X_i^y(t)$ true at $t = 0$ are preserved at each jump of any of the Poisson processes involved, and that $|X^x| - |X^y|$ is decreasing over time. □

The previous results make it possible to establish the fluid regime of the system by studying the rescaled process $\overline{X}_N$ defined by

$$\overline{X}_N(t) = \frac{X(Nt)}{N}, \qquad t \geq 0. \tag{10}$$

In the following, $\overline{L}_N$ denotes the rescaled number of particles, that is, $\overline{L}_N(t) = L(Nt)/N$, and $\overline{\chi}_N = \overline{X}_N/\overline{L}_N$ are the corresponding proportions. Note that any fluid limit is discontinuous at $0^+$ (so that strictly speaking, $X$ does not have any fluid limit), because Proposition 5.1 will show that the fluid limit is at $\pi$ at time $0^+$, and Theorem 6.2 will imply that it stays forever proportional to $\pi$.

Corollary 6.1. *Assume $\lambda > \mu$, and let $x : [0, +\infty[ \mapsto \mathbb{R}^n$ be defined by*

$$x(t) = (1 + (\lambda - \mu)t)\pi.$$

*Then, for any sequence $(x_N, N \geq 1)$ with $|x_N| = N$, any $s, t$ such that $0 < s < t$ and any $\varepsilon > 0$,*

$$\lim_{N \to +\infty} \mathbb{P}_{x_N}\Big(\sup_{s \leq u \leq t} \|\overline{X}_N(u) - x(u)\| \geq \varepsilon\Big) = 0.$$

Proof. Since the size of the initial state goes to infinity, Proposition 5.1 shows that for any $\delta > 0$, the event $\{T_\delta \leq t_\delta\}$ occurs with high probability. Since $T_\delta$ is a stopping time, the strong Markov property makes it possible to use $X_{T_\delta}$ as a new initial point, which is as close to equilibrium as desired. Since, moreover, the total number of customers did not significantly evolve in this time interval, this initial point will satisfy the hypotheses of Theorem 6.2, which makes it possible to conclude.

Denote $\Delta_N(s,t)$ the distance of interest,

$$\Delta_N(s,t) = \sup_{s \leq u \leq t} \|\overline{X}_N(u) - x(u)\|. \tag{11}$$

First, the following decomposition makes it possible to consider all further convergences on the set $\{T_\delta \leq t_\delta\}$:

$$\mathbb{P}_{x_N}(\Delta_N(s,t) \geq \varepsilon) \leq \mathbb{P}_{x_N}(\Delta_N(s,t) \geq \varepsilon, T_\delta \leq t_\delta) + \mathbb{P}_{x_N}(T_\delta > t_\delta)$$

and the last term goes to 0 by Proposition 5.1. The strong Markov property used with the stopping time $T_\delta$ then shows that

$$\mathbb{P}_{x_N}(\Delta_N(s,t) \ge \varepsilon, T_\delta \le t_\delta) \le \mathbb{E}_{x_N}[\mathbb{P}_{X(T_\delta)}(\Delta_N(0,t) \ge \varepsilon)].$$

Now, we isolate the event of interest $\{|L(T_\delta) - |x_N|| \le \sqrt{N}\}$ by writing

$$\begin{aligned}&\mathbb{E}_{x_N}[\mathbb{P}_{X(T_\delta)}(\Delta_N(0,t) \ge \varepsilon); |L(T_\delta) - |x_N|| \le \sqrt{N}] \\ &\quad\le \max_{\substack{y \in \mathbb{N}^n : ||y| - |x_N|| \le \sqrt{N} \\ \text{and } \|y/|y| - \pi\| \le \delta}} \mathbb{P}_y(\Delta_N(0,t) \ge \varepsilon);\end{aligned}$$

therefore, if we note $y_N$ the value that realizes this maximum (the set over which the maximum is considered is finite),

$$\mathbb{E}_{x_N}[\mathbb{P}_{X(T_\delta)}(\Delta_N(0,t) \ge \varepsilon)] \le \mathbb{P}_{x_N}(|L(T_\delta) - |x_N|| \ge \sqrt{N}) + \mathbb{P}_{y_N}(\Delta_N(0,t) \ge \varepsilon).$$

The following inequality holds for any time $u \ge 0$ and any initial state (Lemma 5.1):

$$|L(u) - L(0)| \le \mathcal{N}_\lambda(u) + \mathcal{N}_\mu(u) \stackrel{\text{def}}{=} \mathcal{N}_{\lambda+\mu}(u),$$

and yields

$$\begin{aligned}&\mathbb{P}_{x_N}(|L(T_\delta) - |x_N|| \ge \sqrt{N}) \\ &\quad\le \mathbb{P}_{x_N}(\mathcal{N}_{\lambda+\mu}(T_\delta) \ge \sqrt{N}, T_\delta \le t_\delta) + \mathbb{P}_{x_N}(T_\delta > t_\delta) \\ &\quad\le \mathbb{P}(\mathcal{N}_{\lambda+\mu}(t_\delta) \ge \sqrt{N}) + \mathbb{P}_{x_N}(T_\delta > t_\delta).\end{aligned}$$

This last sum vanishes, so that all is left to prove is that as $N \to +\infty$,

$$\mathbb{P}_{y_N}(\Delta_N(0,t) \ge \varepsilon) = \mathbb{P}_{y_N}\left(\sup_{0 \le u \le t} \|\overline{X}_N(u) - x(u)\| \ge \varepsilon\right) \to 0.$$

Note that the initial state $y_N$ is now such that $|y_N|/N$ goes to 1 (because $|x_N| = N$ and $||y_N| - |x_N|| \le \sqrt{N}$), and $H(y_N/|y_N|, \pi)$ is as small as needed to apply Theorem 6.2, since $\|y_N/|y_N| - \pi\| \le \delta$ and $\delta > 0$ is arbitrary small.

The triangular inequality and the definition of $x$ give for any $0 \le u \le t$,

$$\begin{aligned}\|\overline{X}_N(u) - x(u)\| &\le \|\overline{X}_N(u) - \overline{L}_N(u)\pi\| + \|[\overline{L}_N(u) - (1 + (\lambda - \mu)u)]\pi\| \\ &\le \|\overline{\chi}_N(u) - \pi\| \sup_{0 \le u \le t} \overline{L}_N(u) \\ &\quad + \|\pi\| \sup_{0 \le u \le t} |\overline{L}_N(u) - (1 + (\lambda - \mu)u)|,\end{aligned}$$

and so

$$\mathbb{P}_{y_N}(\Delta_N(0,t) \ge \varepsilon)$$

$$\leq \mathbb{P}_{y_N}\Big(\sup_{0\leq u\leq t}\|\overline{\chi}_N(u)-\pi\| \geq \varepsilon/\Big(2\sup_{0\leq u\leq t}\overline{L}_N(u)\Big)\Big) \tag{12}$$
$$+\,\mathbb{P}_{y_N}\Big(\sup_{0\leq u\leq t}|\overline{L}_N(u)-(1+(\lambda-\mu)u)| \geq \varepsilon/(2\|\pi\|)\Big).$$

Under $\mathbb{P}_{y_N}$, a trivial upper bound for $\overline{L}_N(u)$ for $0 \leq u \leq t$ is given by

$$\overline{L}_N(u) \leq \frac{1}{N}(|y_N| + \mathcal{N}_\lambda(Nt)) \stackrel{\text{def}}{=} A_N(t);$$

therefore, for any constant $C > 0$,

$$\mathbb{P}_{y_N}\Big(\sup_{0\leq u\leq t}\|\overline{\chi}_N(u)-\pi\| \geq \varepsilon/\Big(2\sup_{0\leq u\leq t}\overline{L}_N(u)\Big)\Big)$$
$$\leq \mathbb{P}_{y_N}\Big(\sup_{0\leq u\leq t}\|\overline{\chi}_N(u)-\pi\| \geq \varepsilon/(2C)\Big) + \mathbb{P}(A_N(t)\geq C).$$

For any $t \geq 0$, $A_N(t)$ converges almost surely to $1+\lambda t$ as $N$ goes to infinity; therefore, $\mathbb{P}(A_N(t) \geq C)$ goes to 0 for $C = 2(1+\lambda t)$. The other term vanishes as well. Indeed,

$$\mathbb{P}_{y_N}\Big(\sup_{0\leq u\leq t}\|\overline{\chi}_N(u)-\pi\| \geq \varepsilon/(2C)\Big) = \mathbb{P}_{y_N}(T^{\varepsilon/(2C)} \leq Nt),$$

and by Lemma 2.1, there exists some $\varepsilon' > 0$ such that $T^{\varepsilon/(2C)} \geq T_H^{\varepsilon'}$, hence

$$\mathbb{P}_{y_N}(T^{\varepsilon/(2C)} \leq Nt) \leq \mathbb{P}_{y_N}(T_H^{\varepsilon'} \leq Nt) \leq \mathbb{P}_{y_N}(T_H^{\varepsilon'} < +\infty).$$

One can, moreover, assume $\varepsilon' < \varepsilon_0$ without loss of generality. Observe that so far, $\delta$ is arbitrary; it can be chosen small enough, say $\delta \leq \delta_0$, so that using Lemma 2.1, $H(y_N/|y_N|,\pi) \leq \varepsilon'/2$, and $y_N$ thus satisfies the hypotheses of Theorem 6.2, which shows that $\mathbb{P}_{y_N}(T_H^{\varepsilon'} < +\infty)$, and hence the first term in the upper bound of (12), vanishes in the limit $N \to +\infty$.

The second term of (12) is easier to deal with. We reduce the problem to the event $\{\mathcal{T}_0 = +\infty\}$ by using the following upper bound:

$$\mathbb{P}_{y_N}\Big(\sup_{0\leq u\leq t}|\overline{L}(u)-(1+(\lambda-\mu)u)| \geq \varepsilon/(2\|\pi\|)\Big)$$
$$\leq \mathbb{P}_{y_N}(\mathcal{T}_0 < +\infty) \tag{13}$$
$$+\,\mathbb{P}_{y_N}\Big(\sup_{0\leq u\leq t}|\overline{L}(u)-(1+(\lambda-\mu)u)| \geq \varepsilon/(2\|\pi\|), \mathcal{T}_0 = +\infty\Big).$$

The first term $\mathbb{P}_{y_N}(\mathcal{T}_0 < +\infty)$ in the right-hand side of (13) goes to 0 since, by Lemma 2.2, $\mathbb{P}_{y_N}(\mathcal{T}_0 < +\infty) \leq \mathbb{P}_{y_N}(T_H^{\varepsilon'} < +\infty)$ which has just been proved to vanish as $N \to +\infty$.

Because $L(u) = L(0) + \mathcal{N}_\lambda(u) - \mathcal{N}_\mu(u)$ for all $u \geq 0$ on $\{\mathcal{T}_0 = +\infty\}$, we get the following upper bound for the second term:

$$\mathbb{P}\left(\sup_{0\leq u\leq t}\left|\frac{1}{N}(|y_N| + \mathcal{N}_\lambda(Nu) - \mathcal{N}_\mu(Nu)) - (1 + (\lambda - \mu)u)\right| \geq \varepsilon/(2\|\pi\|)\right)$$

and this term goes to 0 thanks to Doob's inequality. The proof is complete. □

**7. Stability of the subcritical regime.** In this section we consider the subcritical regime $\lambda < \mu$, that is, the case when the input rate is smaller than the maximal output rate. As in the previous section, the key ingredients are the short- and long-term "homogenization" property and the $M/M/1$-like behavior of the total number of customers. The next lemma will be useful for establishing the fluid behavior of the system. It gives a control on the stopping time $T_H^\varepsilon$, or equivalently $T^\varepsilon$; with high probability, $T_H^\varepsilon$ is larger than the time needed for a stable $M/M/1$ queue to empty.

LEMMA 7.1. *Assume $\lambda < \mu$. Fix some $a > 0$ and let $(x_N, N \geq 1)$ be any sequence in $\mathbb{N}^n$ such that*

$$\lim_{N\to+\infty}\frac{|x_N|}{N} = a \quad \textit{and} \quad \lim_{N\to+\infty} H(x_N/|x_N|, \pi) = 0.$$

*Then for any $t < a/(\mu - \lambda)$ and any $\varepsilon < \varepsilon_0$ where $\varepsilon_0$ is given by Lemma 2.2,*

$$\lim_{N\to+\infty}\mathbb{P}_{x_N}(T_H^\varepsilon \leq Nt) = 0.$$

PROOF. Denote $H_N = H(x_N/|x_N|, \pi)$, and let $(\ell_N, N \geq 1)$ be a sequence of integers such that $N \gg \ell_N \gg NH_N$ and $\ell_N \gg \log N$ [such a sequence clearly exists, for example, $\ell_N = N\sqrt{H_N} \vee (\log N)^2$]. Proposition 3.2 with $\alpha = 1/N^2$ and $\delta = \varepsilon/2$ gives

$$\mathbb{E}_{x_N}(e^{-\theta T_H^\varepsilon/N^2}; L(T_H^\varepsilon) \geq \ell_N) \leq C_{\varepsilon/2} e^{2n\log N + |x_N|H_N - \varepsilon\ell_N/2}, \tag{14}$$

where the last bound goes to 0 by choice of $\ell_N$. Let now $\tau_N$ be defined by $\tau_N = \inf\{t \geq 0 : L(t) \leq \ell_N\}$. Since $\ell_N$ is an integer and $L$ has jumps $\pm 1$, we have $L(\tau_N) = \ell_N$, and consequently, for any $t > 0$,

$$\begin{aligned}\mathbb{E}_{x_N}(e^{-\theta T_H^\varepsilon/N^2}; L(T_H^\varepsilon) \geq \ell_N) &\geq \mathbb{E}_{x_N}(e^{-\theta T_H^\varepsilon/N^2}; T_H^\varepsilon \leq \tau_N) \\ &\geq e^{-\theta t/N}\mathbb{P}_{x_N}(T_H^\varepsilon \leq \tau_N \wedge Nt).\end{aligned} \tag{15}$$

Inequalities (14) and (15) together imply that $\mathbb{P}_{x_N}(T_H^\varepsilon \leq \tau_N \wedge Nt)$ goes to 0 as $N$ goes to infinity. Since

$$\mathbb{P}_{x_N}(T_H^\varepsilon \leq Nt) \leq \mathbb{P}_{x_N}(T_H^\varepsilon \leq \tau_N \wedge Nt) + \mathbb{P}_{x_N}(\tau_N < Nt),$$

all is left to prove is that $\mathbb{P}_{x_N}(\tau_N < Nt)$ goes to 0 if $t < a/(\mu-\lambda)$. Using the lower bound $L(t) \geq L(0) + \mathcal{N}_\lambda(t) - \mathcal{N}_\mu(t)$ from Proposition 4.1 and the fact that $|x_N| \geq (\mu-\lambda)Nt + \ell_N$ for $N$ large enough if $t < a/(\mu-\lambda)$, we get for such a $t$,

$$\begin{aligned}
&\mathbb{P}_{x_N}(\tau_N < Nt) \\
&\quad \leq \mathbb{P}_{x_N}(\exists s \in [0, Nt] : L(0) + \mathcal{N}_\lambda(s) - \mathcal{N}_\mu(s) \leq \ell_N) \\
&\quad \leq \mathbb{P}_{x_N}\Big(\sup_{0\leq s\leq Nt} (\mathcal{N}_\mu(s) - \mathcal{N}_\lambda(s) - (\mu-\lambda)s) \geq |x_N| - (\mu-\lambda)Nt - \ell_N\Big) \\
&\quad \leq \mathbb{P}_{x_N}\Big(\sup_{0\leq s\leq Nt} (\mathcal{N}_\mu(s) - \mathcal{N}_\lambda(s) - (\mu-\lambda)s)^2 \geq (|x_N| - (\mu-\lambda)Nt - \ell_N)^2\Big).
\end{aligned}$$

Since $(\mathcal{N}_\mu(s) - \mathcal{N}_\lambda(s) - (\mu-\lambda)s, s \geq 0)$ is a martingale, Doob's inequality yields that the last term is in turn upper bounded by

$$\frac{\operatorname{Var}\mathcal{N}_\lambda(Nt) + \operatorname{Var}\mathcal{N}_\mu(Nt)}{(|x_N| - (\mu-\lambda)Nt - \ell_N)^2} \sim \frac{(\lambda+\mu)Nt}{(a-(\mu-\lambda)t)^2N^2} \to 0. \qquad \square$$

The fluid behavior can now be established. Recall that the rescaled process $\overline{X}_N$ is defined by $\overline{X}_N(t) = X(Nt)/N$ for any $t \geq 0$. In what follows, for $u \in \mathbb{R}$, $u^+$ denotes $\max(u, 0)$.

PROPOSITION 7.1. *Let $x : [0, +\infty[ \mapsto \mathbb{R}^n$ be defined by*

$$x(t) = (1 + (\lambda-\mu)t)^+\pi.$$

*Then for all $0 < s < t$ and all $\varepsilon > 0$,*

$$\lim_{N\to+\infty}\Big[\max_{x\in\mathbb{N}^n : |x|=N} \mathbb{P}_x\Big(\sup_{s\leq u\leq t} \|\overline{X}_N(u) - x(u)\| \geq \varepsilon\Big)\Big] = 0.$$

PROOF. Lemma 7.1 makes it possible to study the system for $t < 1/(\mu-\lambda)$. An additional coupling argument, involving larger initial states, is then required to show that fluid limits stay at 0 after that time. For this technical reason, initial states of size equivalent to $aN$ for some $a > 0$ will be considered, and the following more general result will be established: for $a > 0$, let $x_a : [0, +\infty[ \mapsto \mathbb{R}^n$ be defined by

$$x_a(t) = (a + (\lambda-\mu)t)^+\pi.$$

It will be proved that for any $a > 0$, any $s, t$ with $0 < s < t$ and all $\varepsilon > 0$,

$$\lim_{N\to+\infty}\Big[\max_{x\in\mathbb{N}^n : |x|=\lfloor aN\rfloor} \mathbb{P}_x\Big(\sup_{s\leq u\leq t} \|\overline{X}_N(u) - x_a(u)\| \geq \varepsilon\Big)\Big] = 0,$$

where the notation of the previous section are used.

First assume $t < t_a = a/(\mu - \lambda)$, and set $\Delta_N(s,t) = \sup_{s \le u \le t} \|\overline{X}_N(u) - x_a(u)\|$. The first steps of the proof are similar to the underloaded regime, namely using the strong Markov property to replace the arbitrary initial state by some initial state with low entropy. More precisely, let $\delta_N$ and $t_N$ be as in Proposition 5.2. For any $x \in \mathbb{N}^n$ with $|x| = N$, one has

$$\mathbb{P}_x(\Delta_N(s,t) \ge \varepsilon) \le \mathbb{P}_x(\Delta_N(s,t) \ge \varepsilon, T_{\delta_N} \le Ns) + \mathbb{P}_x(T_{\delta_N} > Ns).$$

Since $t_N/N$ goes to 0, Proposition 5.2 gives that the last term $\mathbb{P}_x(T_{\delta_N} > Ns)$ goes to 0 uniformly in $x \in \mathbb{N}^n$ with $|x| = N$. As for the first term, we write

$$\begin{aligned}
&\mathbb{P}_x(\Delta_N(s,t) \ge \varepsilon, T_{\delta_N} \le Ns) \\
&\quad \le \mathbb{E}_x[\mathbb{P}_{X(T_{\delta_N})}(\Delta_N(0,t) \ge \varepsilon)] \\
&\quad \le \mathbb{P}_{y_N}(\Delta_N(0,t) \ge \varepsilon) + \mathbb{P}_x(|L(T_{\delta_N}) - L(0)| \ge \sqrt{N}),
\end{aligned}$$

where $y_N \in \mathbb{N}^*$ is such that

$$\mathbb{P}_{y_N}(\Delta_N(0,t) \ge \varepsilon) = \max_{\substack{y \in \mathbb{N}^n :\, ||y| - \lfloor aN \rfloor| \le \sqrt{N} \\ \text{and } \|y/|y| - \pi\| \le \delta_N}} \mathbb{P}_y(\Delta_N(0,t) \ge \varepsilon).$$

Because $T_{\delta_N} \le t_N$ with high probability, and because $t_N/\sqrt{N} \to 0$, one can show similarly as in Section 6 that as $N$ goes to infinity,

$$\max_{x \in \mathbb{N}^n :\, |x| = N} \mathbb{P}_x(|L(T_{\delta_N}) - L(0)| \ge \sqrt{N}) \to 0.$$

Along the same lines as in the overloaded case, one gets, by introducing the term $\overline{L}_N(u)\pi$, that for any $C > 0$,

$$\begin{aligned}
&\mathbb{P}_{y_N}(\Delta_N(0,t) \ge \varepsilon) \\
&\quad \le \mathbb{P}_{y_N}\Big(\sup_{0 \le u \le t} \|\overline{\chi}_N(u) - \pi\| \ge \varepsilon/(2C)\Big) \\
&\qquad + \mathbb{P}_{y_N}\Big(\sup_{0 \le u \le t} \overline{L}_N(u) \ge C\Big) \\
&\qquad + \mathbb{P}_{y_N}\Big(\sup_{0 \le u \le t} |\overline{L}_N(u) - (a + (\lambda - \mu)u)| \ge \varepsilon/(2\|\pi\|)\Big).
\end{aligned} \tag{16}$$

Note that since $t < t_a$, $x_a(u) = (a + (\lambda - \mu)u)\pi$ for $0 \le u \le t$. For $C$ large enough, $\mathbb{P}_{y_N}(\sup_{0 \le u \le t} \overline{L}_N(u) \ge C)$ goes to 0 as $N$ goes to infinity. Moreover, Lemma 2.1 gives

$$\mathbb{P}_{y_N}\Big(\sup_{0 \le u \le t} \|\overline{\chi}_N(u) - \pi\| \ge \varepsilon/(2C)\Big) = \mathbb{P}_{y_N}(T^{\varepsilon/(2C)} \le Nt) \le \mathbb{P}_{y_N}(T_H^{\varepsilon'} \le Nt)$$

for some $\varepsilon' > 0$ that can be assumed to satisfy $\varepsilon' < \varepsilon_0$. Since the sequence $(y_N, N \geq 1)$ satisfies the hypotheses of Lemma 7.1 and since $t < t_a$, this last upper bound goes to 0. Moreover, since

$$\begin{aligned}
&\mathbb{P}_{y_N}\Big(\sup_{0\leq u\leq t} |\overline{L}_N(u) - (a + (\lambda-\mu)u)| \geq \varepsilon/(2\|\pi\|)\Big) \\
&\qquad \leq \mathbb{P}_{y_N}(\mathcal{T}_0 \leq Nt) \\
&\qquad\quad + \mathbb{P}_{y_N}\Big(\sup_{0\leq u\leq t} |\overline{L}_N(u) - (a + (\lambda-\mu)u)| \geq \varepsilon/(2\|\pi\|), \mathcal{T}_0 > Nt\Big),
\end{aligned}$$

we conclude, using Lemma 7.1 together with Lemma 2.2 for the first term, and Doob's inequality for the second one, that

$$\mathbb{P}_{y_N}\Big(\sup_{0\leq u\leq t} |\overline{L}_N(u) - (a + (\lambda-\mu)u)| \geq \varepsilon/(2\|\pi\|)\Big) \to 0.$$

The proof in the case $0 < s < t < t_a$ is thus complete.

To conclude in the other cases, a monotonicity argument derived from Lemma 6.2 is used. Let $0 < s < t$ and $t \geq t_a$, and assume in a first step that $t - s < \varepsilon/(2(\mu-\lambda))$. In addition, let $b > (\mu-\lambda)t$ be fixed, and let $t_b = b/(\mu-\lambda)$ be the corresponding time. Note that $t \geq t_a$ implies that $b > a$, so that for any $x \in \mathbb{N}^n$ with $|x| = \lfloor aN \rfloor$, there exists some $y \in \mathbb{N}^n$ such that $y \geq x$ componentwise and $|y| = \lfloor bN \rfloor$. For such $x, y$, Lemma 6.2 shows that $X^x$ and $X^y$ can be coupled in such a way that $|X^x(t)| \leq |X^y(t)|$ for any $t \geq 0$. Hence for any $u \geq s$, using the inequality $\|v\| \leq |v| \leq n\|v\|$ for any $v \in \mathbb{R}^n$, one gets

$$\begin{aligned}
\|\overline{X}^x_N(u) - x_a(u)\| &\leq |\overline{X}^x_N(u)| + |x_a(s)| \leq |\overline{X}^y_N(u)| + |x_a(s)| \\
&\leq n\|\overline{X}^y_N(u) - x_b(u)\| + |x_b(s)| + |x_a(s)|.
\end{aligned}$$

By definition,

$$|x_b(s)| + |x_a(s)| = (\mu-\lambda)(t_b - s) + (\mu-\lambda)(t_a - s)^+ \leq 2(\mu-\lambda)(t_b - s).$$

This yields in turn

$$\mathbb{P}_x\Big(\sup_{s\leq u\leq t} \|\overline{X}_N(u) - x_a(u)\| \geq \varepsilon\Big) \leq \mathbb{P}_y\Big(\sup_{s\leq u\leq t} \|\overline{X}_N(u) - x_b(u)\| \geq \varepsilon''\Big),$$

where $\varepsilon'' = (\varepsilon - 2(\mu-\lambda)(t_b - s))/n$, and finally,

$$\begin{aligned}
&\max_{x\in\mathbb{N}^n : |x|=\lfloor aN\rfloor} \mathbb{P}_x\Big(\sup_{s\leq u\leq t} \|\overline{X}_N(u) - x_a(u)\| \geq \varepsilon\Big) \\
&\qquad \leq \max_{y\in\mathbb{N}^n : |y|=\lfloor bN\rfloor} \mathbb{P}_y\Big(\sup_{s\leq u\leq t} \|\overline{X}_N(u) - x_b(u)\| \geq \varepsilon''\Big).
\end{aligned}$$

Since it has been assumed that $t-s<\varepsilon/(2(\mu-\lambda))$, $b>(\mu-\lambda)t$ can be chosen small enough so that $\varepsilon''>0$. Since $t<t_b$, the first part of the proof implies that

$$\lim_{N\to+\infty}\Big[\max_{y\in\mathbb{N}^n:|y|=\lfloor bN\rfloor}\mathbb{P}_y\Big(\sup_{s\le u\le t}\|\overline{X}_N(u)-x_b(u)\|\ge\varepsilon''\Big)\Big]=0.$$

This proves in particular that when $0<s<t$ and $t-s<\varepsilon/(2(\mu-\lambda))$, then $\max_{|x|=\lfloor aN\rfloor}\mathbb{P}_x(\Delta_N(s,t)\ge\varepsilon)\to0$. It is now left to extend this result to any $s,t$ such that $s<t$, which is a consequence of the following decomposition:

$$\max_{x\in\mathbb{N}^n:|x|=\lfloor aN\rfloor}\mathbb{P}_x(\Delta_N(s,t)\ge\varepsilon)\le\sum_{j=1}^{q}\Big(\max_{x\in\mathbb{N}^n:|x|=\lfloor aN\rfloor}\mathbb{P}_x(\Delta_N(s_{j-1},s_j)\ge\varepsilon)\Big),$$

where $s_0=s<s_1<\cdots<s_q=t$ and $s_j-s_{j-1}<\varepsilon/(2(\mu-\lambda))$ for $1\le j\le q$. Indeed, it has just been shown that each term of this finite sum goes to 0. □

REMARK 7.1. It can be proved that in the critical case $\lambda=\mu>0$, the fluid limit is constant and equal to $\pi$, that is, if $\lambda=\mu>0$, then for all $0<s<t$ and all $\varepsilon>0$,

$$\lim_{N\to+\infty}\Big[\max_{x\in\mathbb{N}^n:|x|=N}\mathbb{P}_x\Big(\sup_{s\le u\le t}\|\overline{X}_N(u)-\pi\|\ge\varepsilon\Big)\Big]=0.$$

This convergence follows readily from Proposition 7.1 and the following coupling. For $0\le\eta\le\lambda$, if $X^\eta$ is a subcritical process with arrival rate $\lambda-\eta$ and departure rate $\lambda=\mu$, then $X$ and $X^\eta$ can be coupled in such a way that $\|\overline{X}(t)-\overline{X}^\eta(t)\|\le N_\eta(t)$ for all $t\ge0$ where $N_\eta$ is a Poisson process with intensity $\eta$.

Note that the behavior of the fluid limit in the critical case does not make it possible to infer the stability or the transience of process $X$. The analogy with the $M/M/1$ queue nevertheless suggests that it could be recurrent null in this case.

In contrast, the behavior of the fluid limit shows that $X$ is ergodic in the subcritical case $\lambda<\mu$.

PROPOSITION 7.2. *When $\lambda<\mu$, the Markov process $X$ is ergodic.*

PROOF. According to Corollary 9.8, page 259 of [12], it is enough to show that for some deterministic time $T>0$,

$$\lim_{N\to+\infty}\max_{x\in\mathbb{N}^n:|x|=N}\mathbb{E}_x(\overline{L}_N(T))=0.$$

Recall that $\overline{L}_N(T) = L(NT)/N$, and let $\varepsilon > 0$ be fixed; then for $x \in \mathbb{N}^n$ with $|x| = N$,

$$\begin{aligned}\mathbb{E}_x(\overline{L}_N(T)) &\le \varepsilon + \mathbb{E}_x(\overline{L}_N(T); \overline{L}_N(T) > \varepsilon)\\ &\le \varepsilon + (1+\lambda T)\mathbb{P}_x(\overline{L}_N(T) > \varepsilon)\\ &\quad + \frac{1}{N}\mathbb{E}_x(\mathcal{N}_\lambda(NT) - \lambda NT; \overline{L}_N(T) > \varepsilon),\end{aligned}$$

where the second inequality comes from $\overline{L}_N(T) \le \overline{L}_N(0) + \mathcal{N}_\lambda(NT)/N$. For any $T \ge 0$, using Cauchy–Schwarz inequality, an upper bound on the last term is given by

$$\frac{1}{N}\mathbb{E}_x(\mathcal{N}_\lambda(NT) - \lambda NT; \overline{L}_N(T) > \varepsilon) \le \frac{1}{N}\mathbb{E}(|\mathcal{N}_\lambda(NT) - \lambda NT|) \le \sqrt{\frac{\lambda T}{N}},$$

so that finally

$$\max_{x \in \mathbb{N}^n : |x| = N} \mathbb{E}_x(\overline{L}_N(T)) \le \varepsilon + (1+\lambda T) \max_{x \in \mathbb{N}^n : |x| = N} \mathbb{P}_x(\overline{L}_N(T) > \varepsilon) + \sqrt{\frac{\lambda T}{N}},$$

and all that is left to prove is that for some $T > 0$, $\max_{|x|=N} \mathbb{P}_x(\overline{L}_N(T) > \varepsilon)$ goes to 0 as $N$ grows to infinity; this is a direct consequence of Proposition 7.1 with $T = 1/(\mu - \lambda)$ since $x(T) = 0$. The proof is now complete. □

## APPENDIX: MARTINGALE CONSTRUCTION

This appendix is devoted to proving Theorem 3.1 which states the existence of a fundamental family of local martingales. In Proposition A.1, we first establish the harmonicity of a special function $g$ which has an integral form. Then a change of variables leads to the local martingale introduced in Theorem 3.1.

**A.1. An integral harmonic function.** The starting point is the generator $\Omega$ of the Markov process $X$ given, for any $x \in \mathbb{N}^n$ and any function $f : \mathbb{R}^n \mapsto \mathbb{R}$, by

$$\begin{aligned}\Omega(f)(x) = {}& \sum_{i=1}^n \lambda_i(f(x+e_i) - f(x)) + \sum_{i=1}^n \mu_i(f(x-e_i) - f(x))\mathbb{1}_{\{x_i > 0\}}\\ &+ \sum_{1 \le i \ne j \le n} q_{ij} x_i (f(x + e_j - e_i) - f(x)).\end{aligned}$$

In addition to the irreducibility of $Q = (q_{ij})_{1 \le i,j \le n}$, we will require that $Q$ is diagonalizable in $\mathbb{C}$, that is, that there exists a set $(\omega_j, 1 \le j \le n)$ of eigenvectors of $Q$ that generate $\mathbb{R}^n$. The complex square matrix $\omega = (\omega_{i,j})_{1\le i,j\le n}$ where $\omega_j = (\omega_{i,j})_{1\le i\le n}$ is invertible.

We can assume without loss of generality that $\omega_n = \mathbb{1}$, denoting by $\mathbb{1}$ the vector in $\mathbb{R}^n$ with all coordinates equal to 1, so that $\omega_n$ is associated to the null eigenvalue; more generally, for $1 \le j \le n$, $\theta_j$ will denote the (possibly complex) eigenvalue associated to $\omega_j$. The negative trace of $Q$ is then given by $-\theta = \sum_1^n \theta_i$ with $\theta > 0$.

In the sequel, $\mathcal{H}$ will denote the hyperplane of $\mathbb{R}^n$ defined by

$$\mathcal{H} = \left\{ v \in \mathbb{R}^n : \sum_{i=1}^n \pi_i v_i = 0 \right\}.$$

For $j = 1, \dots, n-1$, $\omega_j \in \mathcal{H}$ since $Q\omega_j = \theta_j \omega_j$ for $\theta_j \neq 0$ implies (in a matricial form where $\pi$ is a row and $\omega_j$ a column) that $\pi\omega_j = (\theta_j)^{-1}\pi Q\omega_j$ which is 0 since $\pi Q = 0$. These $n-1$ eigenvectors then generate $\mathcal{H}$.

We recall some notation and results of Section 3. $(P_t, t \in \mathbb{R})$ denotes the $Q$-generated Markov semi-group of linear operators in $\mathbb{R}^n : P_t = e^{tQ}$, extended to all real indices $t$ into a group. Each $P_t$ has eigenvalues $e^{\theta_j t}$ and eigenvectors $\omega_j$, $j = 1, \dots, n$. For any $v \in \mathbb{R}^n$ and $t \ge 0$, we define

$$\phi(v,t) = (\phi_i(v,t), 1 \le i \le n) = P_{-t} v.$$

If $v \in \mathbb{R}^n$ and $\varphi(v,\cdot)$ is any primitive of $\sum_{i=1}^n [\mu_i \phi_i(v,\cdot)/(1+\phi_i(v,\cdot)) - \lambda_i \phi_i(v,\cdot)]$ on some open subset $V$ of $\{t \ge 0 : \forall i = 1, \dots, n, 1 + \phi_i(v,t) \neq 0\}$, then the function $h_v(t,x)$ defined by

$$h_v(t,x) = e^{\varphi(v,t)} \prod_{i=1}^n (1 + \phi_i(v,t))^{x_i} \tag{17}$$

is space–time harmonic with respect to $\Omega$ in the domain $V \times \mathbb{N}^{*n}$ (see Proposition 3.1).

The suitable domain of integration for constructing our martingale will be

$$\mathcal{D}(t) = \{ v \in \mathcal{H} : \mathbb{1} + \phi(v,t) > 0 \}, \qquad t \in \mathbb{R},$$

where, for any $u \in \mathbb{R}^n$, $u \ge 0$ (resp. $u > 0$) means that $u_i \ge 0$ (resp. $u_i > 0$) for every $i = 1, \dots, n$.

For each $t \in \mathbb{R}$, $\mathcal{D}(t)$ is an open subset of $\mathcal{H}$. Moreover, it is clear from the definition of $\mathcal{D}(t)$ and from the invariance of $\mathcal{H}$ under the group of operators $(P_s, s \in \mathbb{R})$ that for any $v \in \mathbb{R}^n$ and any $t \ge 0$,

$$v \in \mathcal{D}(t) \quad \Longleftrightarrow \quad P_{-t} v \in \mathcal{D}(0).$$

So for any $t \in \mathbb{R}$, $\mathcal{D}(t) = P_t(\mathcal{D}(0))$. Then since $\mathcal{D}(0) = \{ v \in \mathbb{R}^n : \sum_1^n \pi_i v_i = 0, \mathbb{1} + v > 0 \}$ is clearly bounded, each $\mathcal{D}(t) = P_t(\mathcal{D}(0))$ for $t \in \mathbb{R}$ is bounded as well.

Define the subset $\mathcal{A}$ of $\mathcal{H} \times \mathbb{R}$ by

$$\mathcal{A} = \{(v,t) : t \in \mathbb{R} \text{ and } v \in \mathcal{D}(t)\}.$$

The first step is to show that the following choice of $\varphi$ makes sense:

$$\varphi(v,t) = \int_{-\infty}^{t} \sum_{i=1}^{n} \left( \mu_i \frac{\phi_i(v,s)}{1+\phi_i(v,s)} - \lambda_i \phi_i(v,s) \right) ds \tag{18}$$

for $(v,t) \in \mathcal{A}$. This is the object of the following two lemmas, which will also give some regularity properties of $\varphi$ in view of Proposition A.1.

LEMMA A.1. *If $(v,t) \in \mathbb{R}^n \times \mathbb{R}$ satisfies $\mathbb{1} + \phi(v,t) \geq 0$ and $v \neq -\mathbb{1}$, then $\mathbb{1} + \phi(v,s) > 0$ for all $s < t$. As a consequence,*

$$t > s \quad \Longrightarrow \quad \overline{\mathcal{D}(t)} \subset \mathcal{D}(s), \qquad s,t \in \mathbb{R}, \tag{19}$$

*and*

$$\mathcal{D}(t_0) = \bigcup_{t > t_0} \mathcal{D}(t), \qquad t_0 \in \mathbb{R}. \tag{20}$$

PROOF. Let us first remark that the irreducibility of $Q$ implies that for any $r > 0$ and $(i,j) \in \{1,\dots,n\}$, the probability $P_r(i,j)$ that a Markov process with generator $Q$ initiated at $i$ is in state $j$ at time $r$ is positive. Indeed, if $i = i_0, i_1, \dots, i_k = j$ is a path from $i$ to $j$ such that $q_{i_{l-1},i_l} > 0$ for $l = 1,\dots,k$, then there is a positive probability that the process has exactly followed this path by time $r$.

This implies that $P_r u > 0$ for any $r > 0$ and $u \in \mathbb{R}^n$ such that $u \geq 0$ and $u \neq 0$.

Now let $(v,t)$ satisfy the hypotheses in the lemma, then $\mathbb{1} + \phi(v,t) \neq 0$ since

$$0 \neq \mathbb{1} + v = P_t(\mathbb{1} + P_{-t} v) = P_t(\mathbb{1} + \phi(v,t))$$

and the previous property applied to $u = \mathbb{1} + \phi(v,t)$ and $r = t - s$ for $s < t$ gives

$$\mathbb{1} + \phi(v,s) = \mathbb{1} + P_{-s} v = P_{t-s}(\mathbb{1} + P_{-t} v) = P_{t-s}(\mathbb{1} + \phi(v,t)) > 0.$$

The implication $t > s \Longrightarrow \overline{\mathcal{D}(t)} \subset \mathcal{D}(s)$ follows, noticing that $-\mathbb{1} \notin \mathcal{H}$.

To show (20) for $t_0 \in \mathbb{R}$, note that $\mathcal{D}(t_0)$ contains the right-hand side union by (19), and that the reverse holds since, for $v \in \mathcal{D}(t_0)$, the inequality $\mathbb{1} + \phi(v,t_0) > 0$ extends to some neighborhood of $t_0$. $\square$

LEMMA A.2. (i) *For any* $i \in \{1, \dots, n\}$, *the two integrals*

$$\int_{-\infty}^{0} \phi_i(v,s)\,ds \quad \text{and} \quad \int_{-\infty}^{0} \frac{\phi_i(v,s)}{1+\phi_i(v,s)}\,ds$$

*are well defined for* $v \in \mathcal{D}(0)$, *continuous as functions of* $v$ *on this domain and, respectively, bounded and bounded above on* $\mathcal{D}(0)$.

*The function* $\varphi_0$ *can then be defined on* $\mathcal{D}(0)$ *by*

$$\varphi_0(v) = \int_{-\infty}^{0} \sum_{i=1}^{n} \left( \mu_i \frac{\phi_i(v,s)}{1+\phi_i(v,s)} - \lambda_i \phi_i(v,s) \right) ds$$

*and is continuous and bounded above on* $\mathcal{D}(0)$.

(ii) *The function* $\varphi$ *given by (18) is well defined for* $(v,t) \in \mathcal{A}$ *and satisfies*

$$\varphi(v,t) = \varphi_0(P_{-t}v) \qquad (v,t) \in \mathcal{A}.$$

$\varphi$ *is bounded above on* $\mathcal{A}$ *and continuous with respect to* $v \in \mathcal{D}(t)$ *for fixed* $t \in \mathbb{R}$.

PROOF. (i) Notice that, for fixed $v \in \mathbb{R}^n$, the map $s \mapsto \phi(v,s) = e^{-sQ}v$ is continuous on $\mathbb{R}$ (with values in $\mathbb{R}^n$). Moreover, if $v \in \mathcal{H}$, it has a fast decay as $s$ tends to $-\infty$ as a consequence of the exponential fast convergence of $P_t(i,\cdot)$ to $\pi$ (already used in Section 5).

There exist some positive constants $\eta$ and $C_1$ such that, for any $s \le 0$,

$$\max_{1\le i,j\le n} |P_{-s}(i,j) - \pi_j| \le C_1 \cdot e^{\eta s}. \tag{21}$$

This gives, for $s \le 0$ and $v \in \mathcal{H}$,

$$\|\phi(v,s)\| \le C_2 \cdot e^{\eta s}\|v\|, \tag{22}$$

where $C_2 = nC_1$ which ensures the existence of the vectorial integral $\int_{-\infty}^{0} \phi(v, s)\,ds$ for any $v \in \mathcal{H}$. This integral is continuous with respect to $v$ in $\mathcal{H}$ since, for $v \in \mathcal{H}$,

$$\int_{-\infty}^{0} P_{-s}v\,ds = \int_{-\infty}^{0} (P_{-s} - \Pi)v\,ds = \left( \int_{-\infty}^{0} (P_{-s} - \Pi)\,ds \right) v,$$

where $\Pi$ is the square matrix with all lines equal to $\pi = (\pi_1, \dots, \pi_n)$, and the last matricial integral has a coefficientwise meaning [and is well defined due to (21)]. This shows the integral $\int_{-\infty}^{0} \phi(v,s)\,ds$ as a linear function of $v \in \mathcal{H}$, thus proving its continuity with respect to $v \in \mathcal{H}$. The boundedness of this function on $\mathcal{D}(0)$ follows since $\mathcal{D}(0)$ has compact closure in $\mathcal{H}$.

For the second integral, Lemma A.1 together with the condition $v \in \mathcal{D}(0)$ first ensure that $\mathbb{1} + \phi(v,s) > 0$ for $s \le 0$. The existence of this integral then again follows from the continuity of $s \mapsto \phi(v,s)$ and from the exponential decay in (22).

Let us now begin by proving that it is bounded above on $\mathcal{D}(0)$, writing

$$\int_{-\infty}^{0} \frac{\phi_i(v,s)}{1+\phi_i(v,s)}\,ds = \int_{-\infty}^{-1} \frac{\phi_i(v,s)}{1+\phi_i(v,s)}\,ds + \int_{-1}^{0} \frac{\phi_i(v,s)}{1+\phi_i(v,s)}\,ds$$

and upperbounding each term.

It is easy for the second one, since $1+\phi_i(v,s)>0$ implies that $\phi_i(v,s)/(1+\phi_i(v,s)) \le 1$. In particular,

$$\int_{-1}^{0} \frac{\phi_i(v,s)}{1+\phi_i(v,s)}\,ds \le 1.$$

The first term can be extended to $v \in \overline{\mathcal{D}(0)}$ [again by Lemma A.1 and by the exponential decay in (22)] and can be shown to be bounded on $\overline{\mathcal{D}(0)}$. Indeed, for $v \in \overline{\mathcal{D}(0)}$ and $s \le -1$, $1+\phi_i(v,s)$ is positive and tends to 1 as $s$ tends to $-\infty$ uniformly in $v \in \overline{\mathcal{D}(0)}$ since, by (22),

$$\sup_{v\in\overline{\mathcal{D}(0)}} \|\phi(v,s)\| \le C_2 \cdot e^{\eta s} \sup_{v\in\overline{\mathcal{D}(0)}} \|v\|. \tag{23}$$

Then $1+\phi_i(v,s)$ is bounded below by some positive $\delta$ for $(v,s) \in \overline{\mathcal{D}(0)} \times \,]-\infty, -1]$ (by (23), this is the case on $\overline{\mathcal{D}(0)} \times \,]-\infty,-\kappa]$ for $\kappa$ large enough, and $\overline{\mathcal{D}(0)} \times [-\kappa,-1]$ is compact), and the following bound holds for $v \in \overline{\mathcal{D}(0)}$, using (23):

$$\left|\int_{-\infty}^{-1} \frac{\phi_i(v,s)}{1+\phi_i(v,s)}\,ds\right| \le \frac{C_2}{\delta} \cdot \sup_{u\in\overline{\mathcal{D}(0)}} \|u\| \int_{-\infty}^{-1} e^{\eta s}\,ds = \frac{C_2 e^{-\eta}}{\delta\eta} \sup_{u\in\overline{\mathcal{D}(0)}} \|u\|.$$

Let us now show the continuity of $\int_{-\infty}^{0} \phi_i(v,s)/(1+\phi_i(v,s))\,ds$ with respect to $v \in \mathcal{D}(0)$, using the continuity of $\phi_i(v,s)$ for fixed $s$, together with Lebesgue's dominated convergence theorem. The difficulty is that $\phi_i(v,s)/(1+\phi_i(v,s))$ is not clearly dominated uniformly in $v \in \mathcal{D}(0)$ by some integrable function of $s$ on $]-\infty,0]$, since for $s$ close to 0 and $v$ close to the the portion of $\partial\mathcal{D}(0)$ where $1+v_i=0$, the ratio $\phi_i(v,s)/(1+\phi_i(v,s))$ goes to infinity and is not easily controlled.

It is, however, possible to show local domination, using (20) in the particular case $t_0=0$ and dominating the integrand on each $\mathcal{D}(t), t>0$. This will prove continuity on each $\mathcal{D}(t)$, $t>0$, hence continuity on $\mathcal{D}(0)$ since the $\mathcal{D}(t)$, $t>0$, are open subsets of $\mathcal{D}(0)$. The domination uses the same argument as in the last point; if $t>0$, then

$$\delta(t) = \inf\{1+\phi_i(v,s), (v,s) \in \overline{\mathcal{D}(t)} \times \,]-\infty,0]\}$$

is positive. Then, for any $v \in \mathcal{D}(t)$ and $s \leq 0$,

$$\left| \frac{\phi_i(v,s)}{1+\phi_i(v,s)} \right| \leq \frac{C_2}{\delta(t)} \cdot e^{\eta s} \sup_{u \in \overline{\mathcal{D}(t)}} \|u\|,$$

where the right-hand side is integrable on $]-\infty, 0]$, and hence provides the required domination.

(ii) We use the group structure of $(P_s, s \in \mathbb{R})$ to rewrite both integrals as

$$\int_{-\infty}^{t} \phi_i(v,s)\,ds = \int_{-\infty}^{0} \phi_i(v,s+t)\,ds = \int_{-\infty}^{0} \phi_i(P_{-t}v,s)\,ds$$

and

$$\int_{-\infty}^{t} \frac{\phi_i(v,s)}{1+\phi_i(v,s)}\,ds = \int_{-\infty}^{0} \frac{\phi_i(v,s+t)}{1+\phi_i(v,s+t)}\,ds = \int_{-\infty}^{0} \frac{\phi_i(P_{-t}v,s)}{1+\phi_i(P_{-t}v,s)}\,ds,$$

which ensures their existence for $P_t v \in \mathcal{D}(0)$, that is, $v \in \mathcal{D}(t)$ or equivalently $(v,t) \in \mathcal{A}$, and proves the connexion between $\varphi$ and $\varphi_0$, hence the stated properties of $\varphi$. $\square$

The function $\varphi$ is now legitimately defined by (18) for $t \in \mathbb{R}$ and $v \in \mathcal{D}(t)$.

For fixed $t \in \mathbb{R}$, $\varphi(\cdot, t)$ is continuous (hence measurable) with respect to $v$ in $\mathcal{D}(t)$ and bounded above on this domain.

Reversely, for fixed $v \in \mathcal{D}(0)$, $\varphi(v, \cdot)$ is clearly $\mathcal{C}^1$ on the interval $]-\infty, t_v[$, where $t_v = \sup\{t \geq 0 : \mathbb{1} + \phi(v,t) > 0\} > 0$ [by continuity of $\phi(v,\cdot)$ on $\mathbb{R}$], and $\partial\varphi(v,t)/\partial t = \sum_1^n [\mu_i \phi_i(v,t)/(1+\phi_i(v,t)) - \lambda_i \phi_i(v,t)]$ on this interval.

The following proposition constitutes the first step in defining a new space–time harmonic function, obtained by integrating with respect to $v \in \mathcal{D}(t)$ the parametrized family of functions $h_v$ given by (17), with $\varphi$ given by (18).

As $\mathcal{D}(t) \subset \mathcal{H}$ and $\mathcal{H}$ is an $(n-1)$-dimensional subspace of $\mathbb{R}^n$ which is isomorphic to $\mathbb{R}^{n-1}$ through the one-to-one linear mapping

$$H: \begin{array}{ccc} \mathbb{R}^{n-1} & \longrightarrow & \mathcal{H} \\ u & \longmapsto & \hat{u} = \left(u, -\sum_{i=1}^{n-1} \pi_i u_i / \pi_n\right) \end{array}$$

the new harmonic function will rather take the form of an integral over the following subset of $\mathbb{R}^{n-1}$:

$$\mathcal{C}(t) = H^{-1}(\mathcal{D}(t)) = \{u \in \mathbb{R}^{n-1} : \hat{u} \in \mathcal{D}(t)\} = \{u \in \mathbb{R}^{n-1} : \mathbb{1} + \phi(\hat{u}, t) > 0\}.$$

PROPOSITION A.1. *For any locally Lebesgue-integrable $f$ on $\mathbb{R}^{n-1}$, the function $g(t,x)$ given by the formula*

$$\begin{aligned} g(t,x) &= \int_{\mathcal{C}(t)} h_{\hat{u}}(t,x) f(u)\,du \\ &= \int_{\mathcal{C}(t)} e^{\varphi(\hat{u},t)} \cdot \prod_{i=1}^{n} (1+\phi_i(\hat{u},t))^{x_i} \cdot f(u)\,du \end{aligned} \tag{24}$$

*is space–time harmonic in the domain* $[0,+\infty[\times\mathbb{N}^{*n}$.

PROOF. $g$ is well defined on $[0,+\infty[\times\mathbb{N}^{*n}$. Indeed $\prod_1^n(1+\phi_i(\hat{\cdot},t))^{x_i}$ is continuous on $\mathbb{R}^{n-1}$ for fixed $t\geq 0$ and $x\in\mathbb{N}^{*n}$, hence bounded on the bounded set $\mathcal{C}(t)$ [since $\mathcal{C}(t)$ corresponds to $\mathcal{D}(t)$ through $H^{-1}$]. By Lemma A.2, $e^{\varphi(\hat{\cdot},t)}$ also is continuous and bounded on $\mathcal{C}(t)$ since $\varphi(\cdot,t)$ is continuous and bounded above on $\mathcal{D}(t)$. Then since $f$ is locally integrable, the product of these three functions is integrable on $\mathcal{C}(t)$.

We have to show that $\partial g/\partial t$ exists and satisfies

$$\partial g(t,x)/\partial t + \Omega(g(t,\cdot))(x) = 0.$$

As a rough argument, one expects that

$$\frac{\partial g}{\partial t}(t,x) = \int_{\mathcal{C}(t)} \frac{\partial h_{\hat{u}}}{\partial t}(t,x)\cdot f(u)\,du. \tag{25}$$

Indeed the additional derivation term resulting from the $t$-dependency of the domain $\mathcal{C}(t)$ is bound to vanish, since $h_{\hat{u}}(t,x)$ is zero for $u$ on the frontier of $\mathcal{C}(t)$ (recall that $x\in\mathbb{N}^{*n}$). Therefore, since $\Omega$ commutes with integration,

$$\frac{\partial g}{\partial t}(t,x) + \Omega(g(t,\cdot))(x) = \int_{\mathcal{C}(t)} \left[\frac{\partial h_{\hat{u}}}{\partial t}(t,x) + \Omega(h_{\hat{u}}(t,\cdot))(x)\right]\cdot f(u)\,du = 0$$

by harmonicity of the functions $h_v$, using Proposition 3.1 with $V=[0,t_v[$.

To make this rigorous, all that is needed is to prove (25), by fixing some arbitrary $x\in\mathbb{N}^{*n}$ and $t_0\geq 0$, and studying the ratio $[g(t_0+\delta,x)-g(t_0,x)]/\delta$ as $\delta$ tends to zero. The monotonicity of the family of sets $\mathcal{C}(t)$ forces to distinguish the two cases $\delta>0$ and (for $t_0>0$) $\delta<0$. For the sake of shortness we will only present here the case $\delta>0$, the other side being similar. Both cases make a repeated use of the mean-value theorem and of Lebesgue's dominated convergence theorem.

To simplify notation, define for $u\in\mathcal{C}(t), t\geq 0$ and $x\in\mathbb{N}^{*n}$,

$$h(u,t,x) = h_{\hat{u}}(t,x) = e^{\varphi(\hat{u},t)}\prod_{i=1}^{n}(1+\phi_i(\hat{u},t))^{x_i}.$$

Note that $h$ inherits the derivability properties of $\varphi$ with respect to $t$ (the factor involving $\phi$ being $\mathcal{C}^1$ on $\mathbb{R}$), for fixed $u \in \mathcal{C}(0)$ and $x \in \mathbb{N}^{*n}$, $h$ is $\mathcal{C}^1$ on $]-\infty, t_{\hat{u}}[$.

Let us $x \in \mathbb{N}^{*n}$ and $t_0 \geq 0$ be fixed. For any positive $\delta$, using the inclusion $\mathcal{D}(t_0+\delta) \subset \mathcal{D}(t_0)$ one can write,

$$\frac{g(t_0+\delta, x) - g(t_0, x)}{\delta} = \int_{\mathcal{C}(t_0+\delta)} \frac{h(u, t_0+\delta, x) - h(u, t_0, x)}{\delta} f(u)\,du - \int_{\mathcal{C}(t_0)\setminus\mathcal{C}(t_0+\delta)} \frac{h(u, t_0, x)}{\delta} f(u)\,du.$$

Let us first show that the first term tends to $\int_{\mathcal{C}(t_0)} \frac{\partial h}{\partial t}(u, t_0, x) f(u)\,du$ as $\delta$ tends to zero. Using the mean value theorem, since $h(u, \cdot, x)$ is $\mathcal{C}^1$ on $[0, t_0+\delta]$ for $u \in \mathcal{C}(t_0+\delta)$, this first term can be rewritten as

$$\int_{\mathcal{C}(t_0+\delta)} \frac{\partial h}{\partial t}(u, t_0 + p(u)\delta, x) f(u)\,du$$

for some $p(u) \in ]0,1[$ depending on $u, t_0, x$ and $\delta$.

As $\delta$ goes to zero, $\frac{\partial h}{\partial t}(u, t_0+p(u)\delta, x) f(u)$ tends to $\frac{\partial h}{\partial t}(u, t_0, x) f(u)$ and the indicator function of $\mathcal{C}(t_0+\delta)$ tends to the indicator function of $\mathcal{C}(t_0)$ due equation (20) which obviously extends to the sets $\mathcal{C}(t)$. The convergence of the first term will then result from Lebesgue's theorem by computing [we omit the variable $(\hat{u}, t)$ under $\phi$]

$$\frac{\partial h}{\partial t}(u, t, x) = e^{\varphi(\hat{u}, t)} \sum_{i=1}^{n} (1+\phi_i)^{x_i - 1} \prod_{j \neq i} (1+\phi_j)^{x_j} \left( \mu_i - \lambda_i \phi_i (1+\phi_i) + x_i \frac{\partial \phi_i}{\partial t} \right),$$

and then using the following domination: for $0 < \delta < 1$ and $u \in \mathbb{R}^{n-1}$,

$$\left| \frac{\partial h}{\partial t}(u, t_0 + p(u)\delta, x) f(u) \mathbb{1}_{\mathcal{C}(t_0+\delta)} \right| \leq k(t_0, x) \cdot |f(u)| \cdot \mathbb{1}_{\mathcal{C}(t_0)},$$

where the right-hand side is integrable on $\mathbb{R}^{n-1}$, and $k(t_0, x)$ holds for

$$\sup_{\mathcal{A}} e^{\varphi} \times \sup_{\mathcal{D}(t_0)\times[0, t_0+1]} \left| \sum_{i=1}^{n} (1+\phi_i)^{x_i-1} \prod_{j\neq i} (1+\phi_j)^{x_j} \left( \mu_i - \lambda_i \phi_i (1+\phi_i) + x_i \frac{\partial \phi_i}{\partial t} \right) \right|.$$

The convergence of the first term is thus proved.

We now prove that the second term vanishes as $\delta$ tends to 0. For any $u \in \mathcal{C}(t_0) \setminus \mathcal{C}(t_0+\delta)$ there exists some index $i$ (depending on $u$) such that $1+\phi_i(\hat{u}, t_0) > 0$ while $1+\phi_i(\hat{u}, t_0+\delta) \leq 0$, and this implies by the mean value theorem that $0 < 1+\phi_i(\hat{u}, t_0) \leq -\delta \frac{\partial \phi_i}{\partial t}(\hat{u}, t_0 + q(u)\delta)$ for some $q(u) \in ]0,1[$ depending on $u, t_0$ and $\delta$. One can deduce the following upper bound, again assuming $0 < \delta < 1$:

$$\left| \int_{\mathcal{C}(t_0)\setminus\mathcal{C}(t_0+\delta)} \frac{h(u,t_0,x)}{\delta} f(u)\,du \right| \le k \int_{\mathcal{C}(t_0)\setminus\mathcal{C}(t_0+\delta)} |f(u)|\,du,$$

where $k$ is the following constant:

$$\sup_{\mathcal{A}} e^{\varphi} \times \max_{1\le i\le n} \left\{ \sup_{\mathcal{D}(t_0)\times[0,t_0+1]} \left| \frac{\partial \phi_i}{\partial t} \right| \times \sup_{\mathcal{D}(t_0)\times\{t_0\}} \left| (1+\phi_i)^{x_i-1} \prod_{j\ne i} (1+\phi_j)^{x_j} \right| \right\}.$$

The right-hand side of the previous inequality converges to zero, again by Lebesgue's theorem, because $f$ is integrable on the bounded set $\mathcal{C}(t_0)$, and the sets $\mathcal{C}(t_0)\setminus\mathcal{C}(t_0+\delta)$ decrease to $\varnothing$ as $\delta$ decreases to zero, due to relation (20). □

**A.2. Change of variables.** The last step is now a change of variable in the harmonic function given by the integral (24), for a suitable choice of $f$ so as to separate the time and space variables.

It informally consists in choosing as new variables the quantities $\pi_i(1+\phi_i(v,t))$ $(1 \le i \le n)$, changing the domain $\mathcal{D}(t)$ into $\mathcal{P} \equiv \{v \in \mathbb{R}^n : v > 0$ and $\sum_i^n v_i = 1\}$. Formally, it will be slightly more complicated due to integration with respect to Lebesgue's measure on subdomains of $\mathbb{R}^{n-1}$ [the $\mathcal{C}(t)$'s], which forces a round trip from $\mathbb{R}^{n-1}$ through $\mathbb{R}^n$. So to be correct, this change of variable will rather transform the domain $\mathcal{C}(t)$ into the following one:

$$\mathcal{S} = \left\{ u \in \mathbb{R}^{n-1} : u > 0 \text{ and } \sum_{i=1}^{n-1} u_i < 1 \right\}. \tag{26}$$

We need to introduce some additional notation. Denote by $\Delta$ the diagonal $n \times n$ square matrix having $\pi_1, \ldots, \pi_n$ as its diagonal elements, by $J$ the projection

$$J : \begin{array}{ccc} \mathbb{R}^n & \longrightarrow & \mathbb{R}^{n-1} \\ (v_1, \ldots, v_n) & \longmapsto & (v_1, \ldots, v_{n-1}) \end{array}$$

and by $\mathcal{K}$, the hyperplane of $\mathbb{R}^n$ defined by

$$\mathcal{K} = \left\{ v \in \mathbb{R}^n : \sum_{i=1}^{n} v_i = 1 \right\}.$$

$\mathcal{K}$ corresponds to $\mathbb{R}^{n-1}$ through the one-to-one affine transformation (analogous to $H$ from $\mathbb{R}^{n-1}$ to $\mathcal{H}$),

$$K : \begin{array}{ccc} \mathbb{R}^{n-1} & \longrightarrow & \mathcal{K} \\ u & \longmapsto & \tilde{u} = \left( u, 1 - \sum_{i=1}^{n-1} u_i \right). \end{array}$$

Notice that the inverse mapping of $K$ (resp. $H$) is given by the restriction of $J$ to $\mathcal{K}$ (resp. $\mathcal{H}$). The announced change of variable is given by the $t$-depending transformation

$$\Psi_t : \begin{array}{ccc} \mathbb{R}^{n-1} & \longrightarrow & \mathbb{R}^{n-1} \\ u & \longmapsto & J\Delta(P_{-t}Hu + \mathbb{1}). \end{array}$$

The next lemma shows that $\Psi_t$ can be considered for a change of variables.

LEMMA A.3. *For any $t \geq 0$, $\Psi_t$ is a one-to-one affine transformation on $\mathbb{R}^{n-1}$ which inverse mapping is given by*

$$\Psi_t^{-1}(u) = JP_t(\Delta^{-1}Ku - \mathbb{1})$$

*and which Jacobian is* $\mathrm{Jac}(\Psi_t) = e^{\theta t} \prod_{i=1}^{n-1} \pi_i$. *Moreover,* $\Psi_t(\mathcal{C}(t)) = \mathcal{S}$.

PROOF. Since $\Psi_t$ is clearly an affine transformation in $\mathbb{R}^{n-1}$, its Jacobian is the one of its linear part $J\Delta P_{-t}H$. Now $J\Delta = \Delta' J$ where $\Delta'$ is the diagonal $(n-1) \times (n-1)$ square matrix having $\pi_1, \ldots, \pi_{n-1}$ as its diagonal elements, so that

$$\mathrm{Jac}(J\Delta P_{-t}H) = \left(\prod_{i=1}^{n-1} \pi_i\right) \mathrm{Jac}(JP_{-t}H) = \left(\prod_{i=1}^{n-1} \pi_i\right) \mathrm{Jac}(P_{-t}) = e^{\theta t} \prod_{i=1}^{n-1} \pi_i.$$

The second equality results from the facts that $J$ restricted to $\mathcal{H}$ coincides with $H^{-1}$ and that $\mathcal{H}$ is generated by the first $n-1$ eigenvectors of $P_{-t}$, so that $\mathcal{H}$ is invariant under $P_{-t}$, which restriction to $\mathcal{H}$ has Jacobian $\prod_{i=1}^{n-1} e^{-\theta_i t} = \prod_{i=1}^{n} e^{-\theta_i t} = \mathrm{Jac}(P_{-t})$ $(= e^{\theta t})$ since $\theta_n = 0$. In particular $\mathrm{Jac}(\Psi_t) \neq 0$ so that $\Psi_t$ is invertible.

The formula for $\Psi^{-1}$ easily results from the fact that $\Delta(\cdot + \mathbb{1})$ maps $\mathcal{H}$ onto $\mathcal{K}$ and that the inverse mapping of $K$ is given by the restriction of $J$ to $\mathcal{K}$.

Now using $\mathcal{C}(t) = \{u \in \mathbb{R}^{n-1} : \mathbb{1} + P_{-t}Hu > 0\}$ together with the facts that $\Delta$ preserves the relation $v > 0$ and (again) that $\Delta(\cdot + \mathbb{1})$ maps $\mathcal{H}$ onto $\mathcal{K}$, one gets that $\Psi_t(\mathcal{C}(t))$ is included in $J(\{v \in \mathcal{K} : v > 0\}) = \mathcal{S}$. Equality results from a similar argument for $\Psi_t^{-1}(\mathcal{S}) \subset \mathcal{C}(t)$. $\square$

The transformation $\Psi_t$ hence corresponds to the two following diagrams:

$$\begin{array}{ccc} \mathcal{H} \subset \mathbb{R}^n & \xrightarrow{\Delta(P_{-t}\cdot+\mathbb{1})} & \mathcal{K} \subset \mathbb{R}^n \\ H \uparrow & & \downarrow J \\ \mathcal{C}(t) \subset \mathbb{R}^{n-1} & \overset{\Psi_t}{\dashrightarrow} & \mathcal{S} \subset \mathbb{R}^{n-1} \end{array} \qquad \begin{array}{ccc} \mathcal{H} \subset \mathbb{R}^n & \xleftarrow{P_t(\Delta^{-1}\cdot-\mathbb{1})} & \mathcal{K} \subset \mathbb{R}^n \\ J \downarrow & & \uparrow K \\ \mathcal{C}(t) \subset \mathbb{R}^{n-1} & \overset{\Psi_t^{-1}}{\dashleftarrow} & \mathcal{S} \subset \mathbb{R}^{n-1} \end{array}$$

All that is left now is to choose for (24) a family of locally integrable functions in $\mathbb{R}^{n-1}$ which behave nicely with respect to the change of variable $\Psi_t$. It will be given by the functions $f^{\alpha-1}$ for positive $\alpha$'s, where, for $u \in \mathbb{R}^{n-1}$ and $v \in \mathbb{R}^n$,

$$f(u) = \psi(Hu) = \psi(\hat{u}) \quad \text{and} \quad \psi(v) = \prod_{i=1}^{n-1} |(\omega^{-1}v)_i|.$$

Here, for any $v \in \mathbb{R}^n$ and $1 \le i \le n$, $v_i$ denotes the $i$th coordinate of $v$, so that the $(\omega^{-1}v)_i$'s $(1 \le i \le n-1)$ are the first $n-1$ coordinates of $v$ in the base $(\omega_1, \dots, \omega_n)$ of eigenvectors of $Q$. As will become clear in (29), the next lemma establishes the property of $\psi$ that makes it behave nicely with respect to the change of variables given by $\Psi_t$ by isolating the dependency in time in a separate factor:

LEMMA A.4. *For any $t \ge 0$ and $v \in \mathbb{R}^n$, $\psi(P_t v) = e^{-\theta t}\psi(v)$.*

PROOF. This result stems from diagonalizing $P_t$ as $\omega^{-1}P_t\omega = e^{-t\Theta}$ where $\Theta$ is the diagonal $n \times n$ square matrix having $\theta_1, \dots, \theta_n$ as its diagonal elements. This readily gives $\psi(P_t v) = \prod_{i=1}^{n-1}|(e^{-t\Theta}\omega^{-1}v)_i| = e^{-\theta t}\psi(v)$. $\square$

The main technical point is to establish that $f^{\alpha-1}$ is locally integrable; the next lemma provides in addition a useful upper bound.

LEMMA A.5. *$f^{\alpha-1}$ is locally integrable on $\mathbb{R}^{n-1}$ for any $\alpha > 0$. Moreover, for any compact set $T \subset \mathbb{R}^{n-1}$,*

$$\sup_{0<\alpha\le 1}\left(\alpha^n \int_T f(u)^{\alpha-1}\,du\right) < +\infty. \tag{27}$$

PROOF. If $\alpha \ge 1$, $f^{\alpha-1}$ is continuous on $\mathbb{R}^{n-1}$, hence locally integrable. So consider only the case when $0 < \alpha < 1$.

If the matrix $\omega$ has real coefficients (it can be chosen as such when the eigenvalues $\theta_j$ of $Q$ are real, which is, in particular, the case for a reversible $Q$), $f^{\alpha-1}$ is easily shown to be integrable on any compact set $T$ of $\mathbb{R}^{n-1}$ by operating the change of variable

$$\begin{array}{rcl} \mathbb{R}^{n-1} & \longrightarrow & \mathbb{R}^{n-1} \\ u & \longmapsto & ((\omega^{-1}Hu)_i)_{1\le i\le n-1} = J\omega^{-1}Hu, \end{array}$$

which is linear and one-to-one, and transforms $\int_T f(u)^{\alpha-1}\,du$ into the integral over some compact subset of $\mathbb{R}^{n-1}$ of the locally integrable function $\prod_{i=1}^{n-1}|u_i|^{\alpha-1}$ (up to the Jacobian constant factor). Then by considering $A$ large enough so that $J\omega^{-1}H(T) \subset [-A, A]^{n-1}$ and $A \geq 1$, (27) is obtained from the fact that

$$\int_{[-A,A]^{n-1}} \prod_{i=1}^{n-1} |u_i|^{\alpha-1}\, du = (2A^\alpha)^{n-1}\alpha^{-(n-1)} \leq (2A)^{n-1}\alpha^{-n}.$$

This is not directly possible when $\omega$ has nonreal coefficients. In this case we can show that $f^{\alpha-1}$ is upper bounded by $\prod_{i=1}^{n-1}|(Lu)_i|^{\alpha-1}$ for some invertible $(n-1)\times(n-1)$ square matrix $L$ with real coefficients. The change of variable $u \mapsto Lu$ is then possible, showing (27) in this case similarly as before, which implies the local integrability of $f^{\alpha-1}$ for $0<\alpha<1$.

Since $\alpha < 1$, this amounts to lower bounding $f$ by $\prod_{i=1}^{n-1}|(Lu)_i|$.

Call $C$ the complex invertible $(n-1)\times(n-1)$ square matrix associated to the linear mapping $J\omega^{-1}H$ on $\mathbb{R}^{n-1}$, so that $f(u) = \prod_{j=1}^{n-1}|(Cu)_j|$, and write $C = A + iB$ where $A$ and $B$ are real square matrices. For any $p \in [0,1]$, $u \in \mathbb{R}^{n-1}$ and $j \in \{1,\dots,n-1\}$, the following inequalities hold:

$$\begin{aligned}|(Cu)_j| &\geq \max\{|(Au)_j|, |(Bu)_j|\} \geq |p(Au)_j + (1-p)(Bu)_j| \\ &= |((pA + (1-p)B)u)_j|.\end{aligned}$$

All that is left now is to prove the existence of some $p \in [0,1]$ such that $pA + (1-p)B$ is invertible. It is done through considering the degree $n-1$ polynomial with complex variable, $\det(A + zB)$, which is nonzero at $z = i$ (since $\det C \neq 0$), hence not equal to the null polynomial. It then cannot be zero on the whole real interval $[0,1]$, which gives the result. □

PROOF OF THEOREM 3.1. The previous lemma shows that $f^{\alpha-1}$ is a suitable function to plug in (24); since $\varphi(v,t) = \varphi_0(P_{-t}v)$ for $(v,t) \in \mathcal{A}$, rewriting (24) and using the definition of $f$ gives

$$\begin{aligned}g(t,x) &= \int_{\mathcal{C}(t)} e^{\varphi(\hat{u},t)} \prod_{i=1}^{n} (1 + \phi_i(\hat{u},t))^{x_i} f(u)^{\alpha-1}\, du \\ &= \int_{\mathcal{C}(t)} e^{\varphi_0(P_{-t}Hu)} \prod_{i=1}^{n} (\mathbb{1} + P_{-t}Hu)_i^{x_i} \psi(P_t P_{-t}Hu)^{\alpha-1}\, du.\end{aligned}$$

Expressing $P_{-t}Hu$ through $\Psi_t u$ for $u \in \mathbb{R}^{n-1}$, one gets, since $K$ is the inverse of $J$ restricted to $\mathcal{K}$ and $\Delta(P_{-t}Hu + \mathbb{1}) \in \mathcal{K}$,

$$P_{-t}Hu = \Delta^{-1}K\Psi_t u - \mathbb{1}, \qquad u \in \mathcal{C}(t), \tag{28}$$

so that operating the change of variables given by $\Psi_t$ yields by Lemma A.3

$$g(t,x) = e^{\theta t} \prod_{i=1}^{n-1} \pi_i \int_{\mathcal{S}} e^{\varphi_0(\Delta^{-1}\tilde{u}-\mathbb{1})} \prod_{i=1}^{n} (\Delta^{-1}\tilde{u})_i^{x_i} \psi(P_t(\Delta^{-1}\tilde{u} - \mathbb{1}))^{\alpha-1}\, du. \tag{29}$$

Since $\omega_n = \mathbb{1}$, we have $\psi(v - \mathbb{1}) = \psi(v)$ for any $v \in \mathbb{R}^n$, hence $\psi(P_t(\Delta^{-1}\tilde{u} - \mathbb{1})) = \psi(P_t\Delta^{-1}\tilde{u}) = e^{-\theta t}\psi(\Delta^{-1}\tilde{u})$, where the last equality comes from Lemma A.4. The following function $g'$, which only differs from $g$ by a multiplicative factor, is thus again space–time harmonic:

$$g'(t,x) = e^{\theta t}\int_{\mathcal{S}} G(\tilde{u})\prod_{i=1}^{n}(\Delta^{-1}\tilde{u})_i^{x_i}(e^{-\theta t}\psi(\Delta^{-1}\tilde{u}))^{\alpha-1}\,du$$
$$= e^{-\alpha\theta t}\int_{\mathcal{S}} G(\tilde{u})\prod_{i=1}^{n}\left(\frac{\tilde{u}_i}{\pi_i}\right)^{x_i}\psi(\Delta^{-1}\tilde{u})^{\alpha-1}\,du,$$

where we have defined

$$G(v) = e^{\varphi_0(\Delta^{-1}v - \mathbb{1})}, \qquad v \in \mathcal{P}.$$

Hence defining $F$ as

$$F(v) = \psi(\Delta^{-1}v), \qquad v \in \mathbb{R}^n,$$

yields exactly the local martingale of Theorem 3.1. The second expression (2) is easily obtained. All one needs to do to complete the proof of Theorem 3.1 is to check the announced properties of these two functions $F$ and $G$.

First, $G$ is continuous and bounded on $\mathcal{P}$, since, by Lemma A.2, $\varphi_0$ is continuous and bounded above on $\mathcal{D}(0)$ [if $v \in \mathcal{P}$, then $\Delta^{-1}v - \mathbb{1} \in \mathcal{H}$ and $\Delta^{-1}v > 0$, so that $\Delta^{-1}v - \mathbb{1} \in \mathcal{D}(0)$].

Moreover, $F$ is clearly positive and continuous on $\mathbb{R}^n$, and thus bounded on the bounded subset $\mathcal{P}$ of $\mathbb{R}^n$, and so (3) is the only property left to be checked.

Relation (28) and the fact that $\psi(\cdot - \mathbb{1}) = \psi(\cdot)$, together with the definitions of $F$ and $f$, yield that $F(\tilde{u}) = f(\Psi_0^{-1}u)$ according to the following steps:

$$F(\tilde{u}) = \psi(\Delta^{-1}Ku) = \psi(\Delta^{-1}K\Psi_0(\Psi_0^{-1}u) - \mathbb{1}) = \psi(H\Psi_0^{-1}u) = f(\Psi_0^{-1}u).$$

It follows by the change of variables induced by $\Psi_0$ that, for $\alpha \le 1$,

$$\int_{\mathcal{S}} F(\tilde{u})^{\alpha-1}\,du = \int_{\mathcal{S}} f(\Psi_0^{-1}u)^{\alpha-1}\,du = \operatorname{Jac}(\Psi_0)\int_T f(u)^{\alpha-1}\,du,$$

where $T = \overline{\Psi_0^{-1}(\mathcal{S})}$. (3) then follows using Lemma A.5. $\square$

**Acknowledgments.** We are grateful to Philippe Robert who pointed out this problem to us, and suggested substantial improvements both in the form and in the technical content of this paper.

INRIA
Domaine de Voluceau
78153 Le Chesnay
France
E-mail: florian.simatos@inria.fr

Université Paris 7
Case 7012, 175 rue du Chevaleret
75013 Paris
France
E-mail: danielle.tibi@math.jussieu.fr